                          \def\version{22 August, 2007}                            %

\documentclass[11pt, reqno]{amsart}
\usepackage[centertags]{amsmath}
\usepackage{amsthm, a4, latexsym, amssymb}
\usepackage[latin1]{inputenc}
 \usepackage{color}
\def\@rmrk#1#2{\refstepcounter
    {#1}\@ifnextchar[{\@yrmrk{#1}{#2}}{\@xrmrk{#1}{#2}}}

\makeatletter\@addtoreset{equation}{section}\makeatother

\newfont{\bfit}{cmbxti10 scaled 2000}
\newfont{\biggi}{cmr12 scaled 2000}
\newtheorem{step}{STEP}

\newcommand{\bes}{\begin{step}}
\newcommand{\es}{\end{step}}
 \newcommand{\DeltaD}{\Delta\!^{\mbox{\tiny d}}}
 \newcommand{\lambdaD}{\lambda^{\mbox{\tiny d}}}

 \newcommand{\eps}{\varepsilon}
 
 \newcommand{\supp}{{\rm supp}\,}
 \newcommand{\e}{{\rm e}}
 \renewcommand{\i}{{\rm i}}
 
 \newcommand{\dist}{{\rm dist}}

 \newcommand{\R}{\mathbb{R}}
 \newcommand{\Z}{\mathbb{Z}}
 \newcommand{\N}{\mathbb{N}}

 \renewcommand{\d}{\,\mathrm{d}}
 
 \newcommand{\prob}{\mathbb{P}}
 \newcommand{\Prob}{{\rm Prob }}

 \renewcommand{\P}{\mathbb{P}}
 \newcommand{\E}{\mathbb{E}}
 \newcommand{\one}{\1}

 \newcommand{\Fcal}{{\mathcal F}}
 \newcommand{\Hcal}{{\mathcal H}}
 \newcommand{\Lcal}{{\mathcal L}}
 
 \newcommand{\Ccal}{{\mathcal C}}
 
 \newcommand{\heap}[2]{\genfrac{}{}{0pt}{}{#1}{#2}}
 \newcommand{\sfrac}[2]{\mbox{$\frac{#1}{#2}$}}
 
 \newcommand{\ssup}[1] {{\scriptscriptstyle{({#1}})}}
\def\1{{\mathchoice {1\mskip-4mu\mathrm l}      
{1\mskip-4mu\mathrm l}
{1\mskip-4.5mu\mathrm l} {1\mskip-5mu\mathrm l}}}

 \newcommand{\kommentar}[1]{}

\newenvironment{Proof}[1]
{\vskip0.1cm\noindent{\bf #1}{\hspace*{0.3cm}}}{\vspace{0.15cm}}

\renewcommand{\subsection}{\secdef \subsct\sbsect}
\newcommand{\subsct}[2][default]{\refstepcounter{subsection}
\vspace{0.15cm}
{\flushleft\bf \arabic{section}.\arabic{subsection}~\bf #1  }
\nopagebreak\nopagebreak}
\newcommand{\sbsect}[1]{\vspace{0.1cm}\noindent
{\bf #1}\vspace{0.1cm}}

{\nopagebreak {\hfill{$\diamond$}}\\ }

\newtheorem{theorem}{Theorem}[section]
\newtheorem{lemma}[theorem]{Lemma}

\newtheorem{prop}[theorem]{Proposition}

\theoremstyle{definition}

\def\thebibliography#1{\section*{Bibliography}
  \list%
  {\arabic{enumi}.}
    {\settowidth\labelwidth{[#1]}\leftmargin\labelwidth
    \advance\leftmargin\labelsep
    \parsep0pt\itemsep0pt
    \usecounter{enumi}}
    \def\newblock{\hskip .11em plus .33em minus .07em}
    \sloppy                   
    \sfcode`\.=1000\relax}

\def\P{\prob}

\begin{document}
\title[Potential confinement property of the parabolic
Anderson model]{\Large Potential confinement property of the parabolic
Anderson model}
\author[Gabriela Gr\"uninger and Wolfgang K\"onig]{}
\maketitle
\thispagestyle{empty}
\vspace{-0.5cm}

\centerline{\sc By Gabriela Gr\"uninger\footnote{Institut f\"ur Mathematische Statistik,
Fachbereich Mathematik und Informatik, Einsteinstra\ss e 62,
48149 M\"unster, Germany, {\tt grueninger@math.uni-muenster.de}}$^,$\footnotemark[3] and Wolfgang K\"onig\footnote{Universit\"at Leipzig, Mathematisches Institut, Postfach 10 09 20, D-04009 Leipzig, Germany, {\tt koenig@math.uni-leipzig.de}}$^,$\footnote{Partially supported by the DFG Forschergruppe FOR 718 {\it Analysis and Stochastics in Complex Physical Systems}}}
\renewcommand{\thefootnote}{}
\renewcommand{\thefootnote}{1}

\vspace{-1cm}
\centerline{{Universit\"at M\"unster and
Universit\"{a}t Leipzig}}
\vspace{.5cm}
\centerline{\version}

\begin{quote}{\small }{\bf Abstract.} We consider the parabolic Anderson model, the
Cauchy problem for the heat equation with random potential in $\mathbb{Z}^d$. We use
i.i.d.~potentials $\xi\colon \Z^d \rightarrow \R$ in the third universality
class, namely the class of \textit{almost bounded potentials}, in the classification
of van der Hofstad, K\"onig and M\"orters \cite{HKM06}. This class consists of
potentials whose logarithmic moment generating function is regularly varying with
parameter $\gamma=1$, but do not belong to the class of so-called
double-exponentially distributed potentials studied by G\"artner and Molchanov
\cite{GM98}.

In \cite{HKM06} the asymptotics of the expected total mass was identified in terms
of a variational problem that is closely connected to the well-known logarithmic Sobolev
inequality and whose solution, unique up to spatial shifts,
is a perfect parabola. In the present paper
we show that those potentials whose shape (after appropriate vertical shifting and
spatial rescaling) is away from that parabola contribute only negligibly to the
total mass. The topology used is the strong $L^1$-topology on compacts
for the exponentials of the potential. 
In the course of the proof, we show that any sequence of approximate minimisers 
of the above variational formula approaches some spatial shift of the minimiser, the parabola.

\end{quote}

\pagebreak

\section{Introduction and results}

\subsection{The parabolic Anderson model.}\label{sec-PAM}

\noindent We consider the continuous solution $v\colon[0,\infty)\times \Z^d \to
[0,\infty)$
to the Cauchy problem for the heat equation with random coefficients and localised
initial datum,
    \begin{eqnarray}
    \frac{\partial}{\partial t} v(t,z) & = & \DeltaD v(t,z) + \xi(z) v(t,z),
    \qquad \mbox{ for } (t,z)\in(0,\infty)\times \Z^d,\label{PAM}\\
    v(0,z) & = & \one_0(z),\qquad \mbox{ for } z\in\Z^d.\label{PAMinitial}
    \end{eqnarray}
Here $\xi=(\xi(z) \colon z\in \Z^d)$ is an i.i.d.~random potential
with values in $[-\infty,\infty)$, and $\DeltaD$ is the discrete Laplacian,
    $$
    \DeltaD f(z)= \sum_{y \sim z} \bigl[f(y)-f(z)\bigr], \qquad
    \mbox{ for } z\in \Z^d,\, f\colon\Z^d\to\R.
    $$
The parabolic problem~\eqref{PAM} is called the \emph{parabolic Anderson model}. The
operator $\DeltaD+\xi$ appearing on the right is called the \emph{Anderson
Hamiltonian}; its spectral properties are well-studied in mathematical physics.
Equation~\eqref{PAM} describes a random mass transport through a random field of
sinks and sources, corresponding to lattice points $z$ with $\xi(z)<0$,
respectively, $>0$. There is an interpretation in terms of the expected number of
particles at time $t$ in the site $x$ for a branching process with random
space-dependent branching rates. We refer the reader to \cite{GM90}, \cite{M94} and
\cite{CM94} for more background and to \cite{GK04} for a survey on mathematical
results.

The long-time behaviour of the parabolic Anderson problem is well-studied in the
mathematics and
mathematical physics literature because it is an important example of a model
exhibiting an
\emph{intermittency effect}. This means, loosely speaking, that most of the total
mass of the solution,
\begin{equation}\label{Utdef}
U(t)=\sum_{z\in\Z^d} v(t,z), \qquad\mbox{for } t>0,
\end{equation}
is concentrated on a small number of remote islands, called the {\it intermittent
islands}.
A manifestation of intermittency in terms of the moments of $U(t)$ is as follows. For
$0<p<q$, the main contribution to the $q^{\rm th}$ moment of $U(t)$  comes from islands
that contribute only negligibly to the $p^{\rm th}$ moments. Therefore, intermittency
can be defined by the requirement,
\begin{equation}\label{Intermitt}
\limsup_{t\to\infty}\frac {\langle U(t)^p\rangle^{1/p}}{\langle U(t)^q\rangle^{1/q}}=0,
\qquad \mbox{ for $0<p<q$, }
\end{equation}
where $\langle\,\cdot\,\rangle$ denotes expectation with respect to $\xi$. Whenever
$\xi$ is truly
random, the parabolic Anderson model is intermittent in this sense, see
\cite[Theorem~3.2]{GM90}.

We work under the assumption that all positive exponential moments of
$\xi(0)$ are finite and that the upper tails of $\xi(0)$ possess some mild
regularity property. One of the main results of \cite{HKM06} is that four different
universality classes of long-time behaviours of the parabolic Anderson model can be
distinguished: the so-called double-exponential distribution and some degenerate
version of it studied by G\"artner, Molchanov and K\"onig \cite{GM98}, \cite{GKM05},
bounded from above potentials studied by Biskup and K\"onig \cite{BK01}, and
so-called almost bounded potentials studied by van der Hofstad, K\"onig and
M\"orters \cite{HKM06}. 

In the present paper, we only consider the class of \emph{almost
bounded potentials}, which we will recall in Section~\ref{sec-Almostbounded}.
It is our main purpose to determine those shapes of the random potential $\xi$ that contribute most to the
expectation of the total mass, asymptotically as $t\to\infty$. In other words, we
will find a shifted, rescaled version, $\overline\xi_t$, of $\xi$ and an explicit
deterministic function $\widehat \psi\colon \R^d\to\R$ such that the main
contribution to $\langle U(t)\rangle$ comes from the event $\{\overline\xi_t\approx
\widehat\psi\}$, in a sense that will be specified below. This is what we call a
{\it potential confinement property}; it is a specification of the intermittency
phenomenon for the moments of $U(t)$.

\subsection{Almost bounded potentials.}\label{sec-Almostbounded}

\noindent The class of potentials we will be working with is determined by the
following. We need to introduce the logarithmic moment generating function of
$\xi(0)$ given by
\begin{equation}\label{Def_H}
H(t)=\log \bigl\langle \e^{t\xi(0)}\bigr\rangle, \qquad t\in\R.
\end{equation}

{\bf Assumption (HK).} {\it There is a parameter $\rho\in(0,\infty)$ and a
continuous function $\kappa\colon(0,\infty)\rightarrow (0,\infty)$ with
$\lim_{t\to\infty}\kappa(t)/t= 0$ such that, for all $y \ge 0$,}
\begin{equation}\label{Ann_H}
\underset{t \to \infty}{\lim} \frac{H(yt)-yH(t)}{\kappa(t)}=\rho \cdot y \log y. 
\end{equation}

This is class (iii) of \cite{HKM06}, the class of almost bounded potentials. The
convergence in (\ref{Ann_H}) is uniform in $y\in [0,M]$ for any $M>0$. Both $H$ and
$\kappa$ are regulary varying with index $\gamma=1$. According to \cite[Theorem
3.7.3]{BGT87}, \eqref{Ann_H} is satisfied for $\kappa(t)=H(t)-\int_1^t H(s)/s\,\d
s$. If $\xi$ satisfies \eqref{Ann_H}, then $C\xi$ satisfies  \eqref{Ann_H} with
$\rho$ replaced by $C\rho$, for any $C>0$.

Another important object is the function $\alpha\colon (0,\infty)\rightarrow
(0,\infty)$ defined by
\begin{equation}\label{Def_alpha}
\kappa\Big(\frac{t}{\alpha(t)^d}\Big)=\frac{t}{\alpha(t)^{d+2}},\qquad t\gg1. 
\end{equation}
We also will write $\alpha_t$ instead of $\alpha(t)$. Informally, $\alpha(t)$ is the
order of the diameter of the intermittent islands for the moments. That is, the
expected total mass $\langle U(t)\rangle$ is well-approximated by the sub-sum
$\langle\sum_{|x|\leq R\alpha_t} v(t,x)\rangle$ in a certain sense, after the limits
$t\to\infty$ and afterwards $R\to\infty$ are taken.

\begin{lemma}\label{alpha}
The function $\alpha$ is well defined, up to asymptotic equivalence. Furthermore,
$\lim_{t\to\infty}\alpha(t)= \infty$, and $\alpha$ is slowly varying. In particular,
 $\lim_{t\to\infty}t\alpha_t^{-d} =\infty$.
Furthermore, for any $M>0$,
\begin{equation}\label{H_Asympt}
H\Big(\frac{t}{\alpha_t^d}\cdot y\Big)-y\cdot
H\Big(\frac{t}{\alpha_t^d}\Big)=\frac{t}{\alpha_t^{d+2}}\cdot \rho \cdot y \log
y\cdot(1+o(1)) \quad \text{uniformly in}\quad y \in [0,M]. 
\end{equation}
\end{lemma}

\begin{Proof}{Proof.}
All assertions  besides the last one follow directly from \cite[Prop.~1.2]{HKM06}.
The last one follows from (\ref{Def_alpha}) by substituting $t$ with
$t\alpha_t^{-d}$.
\qed
\end{Proof}

\subsection{Asymptotics for the expected total mass.}

\noindent One of the main results of \cite{HKM06}, see Theorem 1.4, is the
description of the asymptotic behavior of the expected total mass of the parabolic
Anderson model for almost bounded potentials:

\begin{theorem}\label{thm-HKM} Assume that the potential distribution satisfies
Assumption (HK). Then there is a number $\chi\in\R$, depending only on the dimension
$d$ and the parameter $\rho$ appearing in Assumption (HK),  such that
\begin{equation}\label{chiintro}
\lim_{t\to\infty}\frac{\alpha(t)^2}t\log\Big(\langle U(t)\rangle
\e^{-H(t\alpha(t)^{-d})\alpha(t)^d}\Big) =-\chi.
\end{equation}
\end{theorem}

The description of $\chi$ is highly interesting and shows a rich structure, some of
which we want to explore in the present paper. The following objects will play a
crucial role in the following. For $\psi\in\mathcal{C}(\mathbb{R}^d)$  define
\begin{equation}\label{Hlambdadef}
{\mathcal L}(\psi)=\frac{\rho}\e\int_{\R^d} \e^{\frac 1\rho \psi(x)}\, \d
x\qquad\mbox{and}\qquad
\lambda(\psi)=\sup_{\heap{g\in H^1(\R^d)}{\|g\|_2 =1}}\Big\{ \bigl\langle\psi,
g^2\bigl \rangle-\bigl\|\nabla g\bigl\|_2^2\Big\},
\end{equation}
where $H^1(\mathbb{R}^d)$ is the usual Sobolev space,  $\nabla$ the usual
(distributional) gradient and $\langle\cdot,\cdot\rangle$ and  $||\cdot ||_2 $ are
the inner product and the norm on $L^2(\mathbb{R}^d)$. Then $\lambda(\psi)$ is the
top of the spectrum of the operator $\Delta+\psi$ in $H^1(\R^d)$. If $\psi$ decays
at infinity sufficiently fast towards $-\infty$, then $\Lcal(\psi)$ is finite and
$\lambda(\psi)$ is the principal $L^2$-eigenvalue of $\Delta+\psi$ in $\R^d$. Now we
can identify $\chi$ explicitly, see \cite[Prop.~1.11]{HKM06}. 

\begin{lemma}\label{lem-chiident} The limit $\chi$ in \eqref{chiintro} is identified as
\begin{equation}\label{chiident}
\chi=\inf_{\psi\in\Ccal(\R^d)\colon\Lcal(\psi)<\infty}\big[\Lcal(\psi)-\lambda(\psi)\big].
\end{equation}
Furthermore, the infimum is uniquely, up to spatial shifts, attained at the parabola
$$
\widehat\psi(x)=\rho+\rho\frac d2\log\frac\rho\pi-\rho^2 |x|^2,\qquad x\in\R^d.
$$
In particular, $\chi=\rho d(1-\frac 12\log \frac\rho \pi)$.
\end{lemma}

\subsection{Heuristic explanation.}\label{sec-Heur}

The content of Theorem~\ref{thm-HKM}, in combination with Lemma~\ref{lem-chiident},
can heuristically be explained in terms of a large-deviation statement as follows.
Introduce the vertically shifted and rescaled version of the potential $\xi$,
\begin{eqnarray}
\xi_t(z)&=&\xi(z)-\frac{\alpha(t)^d}tH\Big(\frac t{\alpha(t)^d}\Big),\qquad
z\in\Z^d,\label{xitdef}\\
\overline\xi_t(x)&=&\alpha(t)^2\xi_t\big(\lfloor \alpha(t)
x\rfloor\big),\qquad\qquad x\in\R^d.\label{xitbardef}
\end{eqnarray}
Then $\overline \xi_t$ is a random step function $\R^d\to\R$. Using a Fourier
expansion with respect to the eigenfunctions of $\DeltaD+\xi$ in large,
$t$-dependent boxes, one can show that the total mass $U(t)$ is asymptotically equal
to $\exp\{t\lambdaD_{t\log t}(\xi)\}$, where $\lambdaD_{t\log t}(V)$ denotes the
principal eigenvalue of the operator $\DeltaD+V$ in the centred box with radius
$t\log t$ with zero boundary condition, for any potential $V\colon \Z^d\to\R$. Some
technical work is done to show that $\lambdaD_{t\log t}(\xi)$ may asymptotically be
replaced by the eigenvalue $\lambdaD_{R\alpha_t}(\xi)$ in the much smaller box of
radius $R\alpha_t$. More precisely, the replacement error is exponential on the
scale $t/\alpha_t^2$, and its rate vanishes if the limit $R\to\infty$ is eventually
taken. Using \eqref{xitbardef} and asymptotic scaling properties of
$\lambdaD_{R\alpha_t}(\cdot)$, we see that
$$
U(t)\e^{-H(t\alpha(t)^{-d})\alpha(t)^d}\approx
\exp\Big\{t\lambdaD_{R\alpha_t}(\xi_t)\Big\}\approx \exp\Big\{\frac
t{\alpha(t)^2}\lambda_R(\overline\xi_t)\Big\};
$$
where $\lambda_R(\psi)$ denotes the principal eigenvalue of $\Delta+\psi$ in the box
$Q_R=[-R,R]^d$ with zero boundary condition; note that the term
$-H(t\alpha(t)^{-d})\alpha(t)^d$ is absorbed in the vertically shifted potential,
$\xi_t$.

Now we take expectations with respect to the potential and find that the expected
total mass is given in terms of an exponential moment of $\lambda_R(\overline\xi_t)$
on the scale $t\alpha_t^{-2}$. The following lemma is one key property of the
shifted and rescaled potential $\overline\xi_t$ and gives to the functional $\Lcal$
defined in \eqref{Hlambdadef} the meaning of a large-deviation rate function. We
introduce the set $\Fcal(Q_R)$ of all measurable functions $\psi\colon Q_R\to\R$
that are bounded from above.

\begin{lemma}[LDP for $\overline\xi_t$]\label{lem-LDP} Fix $R>0$. Then the
restriction of $(\overline\xi_t)_{t>0}$ to $Q_R$ satisfies a large-deviation
principle with speed $t\alpha_t^{-2}$ and rate function 
\begin{equation}\label{LRdef}
\Lcal_R\colon \Fcal(Q_R)\to\R,\qquad \Lcal_R(\psi)=\frac{\rho}\e\int_{Q_R} \e^{\frac
1\rho \psi(x)}\, \d x,
\end{equation}
with respect to the topology that is induced by test integrals against all
nonnegative continuous functions $Q_R\to[0,\infty)$.
\end{lemma}

\begin{Proof}{Sketch of proof.} We identify the limiting cumulant generating function,
$$
\Lambda_R(f)=\lim_{t\to\infty}\frac{\alpha(t)^2}t\log\Big\langle \exp\Big\{\frac
t{\alpha(t)^2}\int_{Q_R} \overline\xi_t(x) f(x)\,\d x\Big\}\Big\rangle,
$$
for any continuous nonnegative $f\colon Q_R\to[0,\infty)$.
Indeed, we shall show that $\Lambda_R(f)$ exists and is equal to $\Hcal_R(f)=\rho
\int_{Q_R} f(x)\log f(x)\,\d x$. Then the well-known G\"artner-Ellis theorem
\cite[Sect.~4.5.3]{DZ98} yields the result, since $\Hcal_R$ is the Legendre
transform of $\Lcal_R$, see also Lemma~\ref{lemchirepr} below.

An explicit calculation using \eqref{xitdef}, Assumption (HK) and \eqref{Def_alpha}
shows that 
$$
\begin{aligned}
\frac{\alpha(t)^2}t\log &\Big\langle \exp\Big\{\frac t{\alpha(t)^2}\int_{Q_R}
\overline\xi_t(x) f(x)\,\d x\Big\}\Big\rangle\\
&=\rho (1+o(1))\int_{Q_R}\d x\, f(x)\log\Big(\int_{\lfloor
x\alpha_t\rfloor/\alpha(t)+Q_{1/\alpha(t)}} f(y)\,\d y\,\alpha(t)^d\Big)\\
&=\rho (1+o(1))\int_{Q_R} f(x)\log f(x)\,\d x.
\end{aligned}
$$
Obviously, this implies that $\Lambda_R(f)$ exists and equals $\Hcal_R(f)$.
\end{Proof}
\qed

We kept this proof short since we are not going to use Lemma~\ref{lem-LDP} in our
proofs. Loosely speaking, this principle says that
\begin{equation}
\lim_{t\to\infty}\frac {\alpha(t)^2}t\log\Prob\big(\overline\xi_t\approx \psi\mbox{
in }Q_R)=-\Lcal_R(\psi),
\end{equation}
for sufficiently regular functions $\psi$. Using this principle in combination with
Varadhan's lemma \cite[Sect.~4.3]{DZ98} and making $R$ very large, we arrive at
$$
\begin{aligned}
\big\langle U(t)\big\rangle \e^{-H(t\alpha(t)^{-d})\alpha(t)^d}&\approx
\Big\langle\exp\Big\{\frac t{\alpha(t)^2}\lambda_R(\overline\xi_t)\Big\}\Big\rangle
\approx\exp\Big\{\frac t{\alpha(t)^2}\sup_{\psi\in
\Fcal(Q_R)}\big[\lambda_R(\psi)-\Lcal_R(\psi)\big]\Big\}\\
&\approx \exp\Big\{-\chi\frac t{\alpha(t)^2}\Big\}.
\end{aligned}
$$

This ends the heuristic derivation of Theorem~\ref{thm-HKM}. Hence, we see that
there is a competition between two forces for large $R$: the potential tries to keep
the value of the eigenvalue $\lambda_R(\overline\xi_t)$ as high as possible, but has
to pay an amount of $\Lcal_R(\overline\xi_t)$ for doing that. The best contribution
comes from potentials $\overline\xi_t$ that make an optimal compromise, i.e.,
optimize the difference of the two contributions. This is precisely what is
expressed in \eqref{chiident}.

\subsection{Our result: potential confinement.}\label{sec-result}

The purpose of the present paper is to give rigorous substance to the heuristics of
Section~\ref{sec-Heur}. We prove that there is a one-to-one correspondence between
near-by minimisers $\psi$ of the variational formula in \eqref{chiident} and the
contribution to the expected total mass coming from events $\{\overline
\xi_t\approx\psi\}$. More precisely, we prove that the contribution to the expected
total mass that comes from potential shapes outside a neighborhood of any shift of
the parabola $\widehat\psi$ is asymptotically negligible with respect to the full
expectation. 
 
Let us first introduce the topology of potentials we are working with. We write
$Q_R=[-R,R]^d$ for the centred cube of sidelength $2R$. Introduce the distance
\begin{equation}\label{metricdef}
\dist(f_1,f_2)=\sum_{r=1}^\infty 2^{-r}
\phi\Big(\int_{Q_r}\big|f_1(x)-f_2(x)\big|\,\d x\Big),\qquad f_1,f_2\in L^1(\R^d),
\end{equation}
where $\phi(s)=\frac s{1+s}$ for $s>0$. Then $\dist$ metrisizes the topology of
$L^1$-convergence
on every compact subset of $\R^d$. For describing general potential realisations, we
enlarge the space of continuous functions to a much larger function set, the set
$\Fcal$ of all measurable functions $\psi\colon\R^d\to\R$ that are bounded from
above. Now we can formulate our main result, a law of large numbers for
$\overline\xi_t$ defined in \eqref{xitdef}--\eqref{xitbardef} towards the set of
minimizers of the formula in \eqref{chiident}.

\begin{theorem}[Potential confinement]\label{thm-main} Suppose that Assumption (HK)
holds. Then
\begin{equation}\label{HSF}
\lim\limits_{t \to \infty} \frac{\bigl\langle U(t)
\1_{\widehat\Gamma_{t,\eps}}(\overline\xi_t) \bigr\rangle}{\bigl\langle U(t)
\bigr\rangle}=0,
\end{equation}
where
\begin{equation}\label{Gammadef}
\widehat\Gamma_{t,\eps}=\bigcap_{M\in(0,\infty)}\bigcap_{x\in Q_{t\log
t}}\Big\{\psi\in\Fcal\colon\dist\big(\e^{\frac 1\rho(\psi(x+\cdot)\wedge M)},
\e^{\frac 1\rho\widehat \psi(\cdot)}\big)>\eps\Big\}.
\end{equation}
\end{theorem}

Theorem~\ref{thm-main} says that the totality of all potential realisations $\xi$
such that every shift of $\e^{\frac 1\rho(\overline \xi_t\wedge M)}$ is, for any
$M>0$, away from the Gaussian density $\e^{\frac 1\rho\widehat\psi}$ by some
positive amount gives a negligible contribution to the expected total mass. It is
sufficient to consider only shifts by amounts $\leq t\log t$ since the mass coming
from farther away contributes negligibly at time $t$ anyway. It will turn out in the
proof that the quotient on the left hand side of \eqref{HSF} decays exponentially on
the scale $t\alpha(t)^{-2}$. The appearance of the parameter $M$ is necessary since
distances between $\e^{\frac 1\rho\overline \xi_t\wedge M}$ and $\e^{\frac
1\rho\overline \xi_t}$ cannot be controlled on that exponential scale.

\subsection{Comments on the proof.}

The proof of Theorem~\ref{thm-main} has a functional analytic side and a
probabilistic side. On one hand, we show that any sequence of functions that
asymptotically minimise $\Lcal-\lambda$ in \eqref{chiident} converge, after an
appropriate spatial translation, to the minimiser $\widehat \psi$ in the topology
used in Theorem~\ref{thm-main}, and on the other hand we derive effective estimates
for the expectation of the total mass on the event that $\overline\xi_t$ is bounded
away from $\widehat \psi$ in the same sense. The main point is that these two
properties have to be proved in the same topology, which is a non-trivial issue.
Note that the topology we work with is much stronger than the one in which we have a
large-deviation principle, see Lemma~\ref{lem-LDP}. In the literature, other
topologies are considered in which the variational formula in \eqref{chiident} has a
related approximation property (see the remarks at the beginning of
Section~\ref{sec-VarForm}); however these topologies turned out to be not suitable for our probabilistic approach.

The analysis part of the proof will be handled in Section~\ref{sec-VarForm} by more
or less standard methods from analysis. The probabilistic part is treated in
Section~\ref{sec-proofthm}. The large-deviations principle of Lemma~\ref{lem-LDP}
can serve as a guidance only since the topology used in that principle is too weak.
Our proof indeed follows another route, which we informally describe now.

Similarly to the heuristics of Section~\ref{sec-Heur}, we have 
$$
\e^{-H(t\alpha(t)^{-d})\alpha(t)^d}\bigl\langle U(t)
\1_{\widehat\Gamma_{t,\eps}}(\overline\xi_t)
\bigr\rangle\approx\Big\langle\exp\Big\{\frac
t{\alpha(t)^2}\lambda_R(\overline\xi_t)\Big\}\1_{\Gamma_{R,\eps}}(\overline\xi_t)\Big\rangle,
$$
where $\Gamma_{R,\eps}$ is some finite-box approximation of
$\widehat\Gamma_{t,\eps}$. Now we add und subtract the term
$t\alpha(t)^{-2}\rho\log(\frac \e\rho\Lcal_R(\overline\xi_t))$ in the exponent. The
difference term is estimated against the variational formula
$$
-\chi_R(\eps)=\sup_{\psi\in\Gamma_{R,\eps}}\Big(\lambda_R(\psi)-\rho\log\Big(\frac
\e\rho\Lcal_R(\psi)\Big)\Big),
$$
such that we have
$$
\e^{-H(t\alpha(t)^{-d})\alpha(t)^d}\bigl\langle U(t)
\1_{\widehat\Gamma_{t,\eps}}(\overline\xi_t) \bigr\rangle
\leq \e^{-t\alpha(t)^{-2} \chi_R(\eps)}\Big\langle \exp\Big\{\frac
t{\alpha(t)^2}\rho\log\Big(\frac
\e\rho\Lcal_R(\overline\xi_t)\Big)\Big\}\Big\rangle.
$$
With the help of the principle in Lemma~\ref{lem-LDP} and Varadhan's lemma one can
convince oneself that the exponential rate (on the scale $t\alpha(t)^{-2}$) of the
last expectation should be equal to
$$
\sup_\psi\Big(\rho\log\Big(\frac
\e\rho\Lcal_R(\psi)\Big)-\Lcal_R(\psi)\Big)=\sup_{l\in(0,\infty)}\Big(\rho\log\frac{\e
l}\rho-l\Big).
$$
(However, the proof of that fact cannot be done with the help of
Lemma~\ref{lem-LDP}, since the functional $\psi\mapsto \rho\log(\frac
\rho\e\Lcal_R(\psi))$ is not bounded and continuous in the topology used in that
lemma.) The right-hand side is easily seen to be zero with unique minimiser
$\Lcal_R(\psi)=l=\rho$. Hence, the only task that is left to do is to prove that
$\liminf_{R\to\infty}\chi_R(\eps)>\chi$. This is indeed true; it relies on the
representation
$$
-\chi=\sup_{\psi}\Big(\lambda_R(\psi)-\rho\log\Big(\frac \e\rho\Lcal_R(\psi)\Big)\Big);
$$
see \cite{HKM06}. From this estimate, Theorem~\ref{thm-main} follows since the
denominator of \eqref{HSF} has the strictly larger exponential rate $-\chi$ by
Theorem~\ref{thm-HKM}.

\subsection{Remarks on the literature.}\label{sec-Lit}

To the best of our knowledge, the only potential confinement property that has been
proved for the parabolic Anderson model in the literature is in \cite{GKM05} for the
universality class of the double-exponential distribution, including its degenerate 
version. That paper works in the
almost-sure setting and proves that the strictly main contribution to the total mass
$U(t)$ comes from islands in which the potential looks like the maximisers of the
relevant variational formula. That formula is the discrete variant of the formula
appearing in the present paper, i.e., for the discrete Laplace operator on $\Z^d$
instead of the continuous one on $\R^d$. 

There is a \lq dual\rq\ confinement property in the parabolic Anderson model, the
confinement of the path of the random walk in the Feynman-Kac formula, see
\eqref{FKformula} below. This property says that the strictly maximal contribution
to the expected total mass $U(t)$ comes from those random walk paths whose shape,
after appropriate rescaling, resembles the minimisers of the \lq dual\rq\ version of
the characteristic variational problem (see Lemma~\ref{lemchirepr} for the dual
representation of $\chi$ in the case handled in the present paper). This property is
proved in $d=2$ by Bolthausen \cite{B94} in an important special case of the
universality class of potentials that are bounded from above: they assume only the
two values $0$ and $-\infty$ in \cite{B94}. A similar result, also in $d=2$, was
independently derived by Sznitman \cite{S91} for the spatially continuous variant
for Brownian motion in a Poisson trap field. The characteristic  variational problem
is in that
  case
$$
\chi=\inf\big\{\|\nabla g\|_2^2+\rho |\supp(g)|\colon g\in
H^1(\R^d),\|g\|_2=1,\supp(g)\mbox{ compact}\big\};
$$
the function $g^2$ plays the role of the normalised rescaled occupation measures of
the walk, respectively of the Brownian motion. The restriction to $d=2$ was removed
by Povel \cite{P99}, after suitable isoperimetric inequalities derived in the
analysis literature had become known.

\section{Proof of Theorem~\ref{thm-main}}\label{sec-proofthm}

In this section, we prove Theorem~\ref{thm-main}. Recall that we suppose that
Assumption (HK) holds, and recall the parameter $\rho\in(0,\infty)$ form that
assumption. Comparing to Theorem~\ref{thm-HKM}, it is easy to see that the following
proposition immediately implies Theorem~\ref{thm-main}. 

\begin{prop}\label{Hauptprop}
For any $\eps \, >0$,
\begin{equation}\label{Haupt1}
\limsup_{t \to \infty} \frac{\alpha_t^2}{t} \log\Bigl(\e^{-\alpha_t^d
H(t/\alpha_t^d)}\Big\langle U(t) \1_{\widehat\Gamma_{t,\eps}}(\overline
\xi_t)\Big\rangle \Bigr)<-\chi.
\end{equation}
\end{prop}

Indeed, Theorem~\ref{thm-HKM} says that the denominator of \eqref{HSF}, after
inserting the factor $\e^{-\alpha_t^d H(t/\alpha_t^d)}$ both in numerator and
denominator, has the exponential rate $-\chi$, while the rate of the numerator is
strictly smaller, according to Proposition~\ref{Hauptprop} (both on the scale
$t\alpha_t^{-2}$). Hence, Proposition~\ref{Hauptprop} implies that the quotient in
\eqref{HSF} even decays exponentially on the scale $t\alpha_t^{-2}$.

One of the most important tools in the study of the parabolic Anderson model is the
{\it Feynman-Kac formula}, which represents the solution of \eqref{PAM} and its
total mass in terms of an exponential expectation of a functional of simple random
walk $(X(s) \colon s\in[0,t])$ on $\Z^d$ with generator
$\DeltaD$. We denote by $\P_z$ and $\E_z$ probability and expectation with respect
to the random walk, when started at $z\in\Z^d$. The walker's {\it local times} are
denoted by $\ell_t(z)=\int_0^t\delta_z(X(s))\, \d s$, the amount of time the walker
spends at $z\in\Z^d$ by time $t>0$. Note that $\int_0^t V(X(s))\d s=\langle
V,\ell_t\rangle$ for functions $V\colon \Z^d\to\R$, where $\langle
f,g\rangle=\sum_{z\in\Z^d}f(z)g(z)$ for any $f$, $g$. Then, also using
\eqref{xitdef}, the Feynman-Kac formula may be formulated by saying 
\begin{equation}\label{FKformula}
\e^{-\alpha_t^d H(t/\alpha_t^d)}U(t)= \E_0\Big[\exp\Big\{\int_0^t \xi_t(X(s))\,\d
s\Big\}\Big]
=\E_0\Big[\e^{\langle \ell_t,\xi_t\rangle}\Big].
\end{equation}

We divide the proof of Proposition~\ref{Hauptprop} into a sequence of steps. In
Section~\ref{sec-boxreduction} we show how we reduce the infinite state space $\Z^d$
to some finite large box. In Section~\ref{sec-potentialcutting} we replace the
shifted and rescaled potential, $\overline\xi_t$, by a truncated version
$\overline\xi_t\wedge M$ and show that the replacement error vanishes as
$M\to\infty$. This technical step turns out to be crucial in
Section~\ref{sec-PropProof} since our proof of Lemma~\ref{Lemma_int} would fail for
$\overline\xi_t$ in place of $\overline\xi_t\wedge M$. After the two preparatory
steps in Sections~\ref{sec-boxreduction} and \ref{sec-potentialcutting}, the main
strategy of the proof of Proposition~\ref{Hauptprop} is carried out in
Section~\ref{sec-PropProof}.

\subsection{Reduction to a large box.}\label{sec-boxreduction}

Our first main step is to estimate the expectation on the left-hand side of
\eqref{Haupt1} in terms of a finite-box version. In other words, we argue that we
may replace the full state space, $\Z^d$, by a box with a radius of order
$\alpha_t$. We will also insert an appropriate scaling, which will turn the discrete
box of order $\alpha_t$ into continuous cubes of finite-order radius. By
$B_R=[-R,R]^d\cap\Z^d$ and $Q_R=[-R,R]^d$ we denote the discrete box and the
continuous cube of radius $\lfloor R\rfloor$, resp.~$R$. A finite-cube version of
the distance $\dist$ defined in \eqref{metricdef}, appropriate for our purposes, is 
\begin{equation}
\d_R(\psi_1,\psi_2)=\int_{Q_{R}}\big|{\rm e}^{\frac 1\rho\psi_1(x)}-{\rm e}^{\frac
1\rho\psi_2(x)}\big|\,\d x,\qquad \psi_1,\psi_2\in \Fcal,
\end{equation}
where we recall that $\Fcal$ denotes the set of all measurable functions $\R^d\to\R$
that are bounded from above.

\begin{lemma}[Reduction to a large box]
\label{Compact}Fix $\eps>0$.
Then there is $C>0$ such that for all $R \ge 2-\log\frac \eps2 $ and $t \gg 1$,
\begin{equation}\label{compact}
\begin{aligned}
{\rm e}^{-\alpha_t^d H(t/\alpha_t^d)}&\Bigl\langle U(t)
\1_{\widehat\Gamma_{t,\eps}}(\overline\xi_t) \Bigr\rangle\\
&\le {\rm e}^{H(2t)/2-\alpha_t^d H(t/\alpha_t^{d})}{\rm e}^{-\frac 12 t\log t}+{\rm
e}^{\frac{C t}{R^2 \alpha_t^2}}\Bigl\langle \E^{t,R}\big[{\rm e}^{\langle \ell_t,
\xi_t\rangle}\big]\1_{\{\forall M>0\colon \overline\xi_t\wedge
M\in\Gamma_{R,\eps/2}\}}\Bigr\rangle,
\end{aligned}
\end{equation}
where we abbreviate $\E^{t,R}[\ldots]=\E_0[\ldots \cdot \1_{\{\supp \ell_t \subset
B_{3R{\alpha_t}}\}}]$, and we put
\begin{equation}\label{GammaRepsdef}
\Gamma_{R,\eps}=\bigcap_{x\in Q_{2R}}\Big\{\psi\in\Fcal\colon
\d_R\big(\psi(x+\cdot), \widehat\psi(\cdot)\big)>\eps\Big\},\qquad \eps>0.
\end{equation}
\end{lemma}

Lemma~\ref{Compact} reduces the $\overline\xi_t$-expectation to an expectation of
the restriction to the cube $Q_{3R}$; the constraint that any shift is away from
$\widehat\psi$ on any compact subset of $\R^d$ is replaced by the requirement that
the shift by any amount $\leq 2R$ is away from the $Q_R$-restriction of
$\widehat\psi$ in $L^1(Q_R)$-sense.

\begin{Proof}{Proof of Lemma~\ref{Compact}.} This is a refinement of the proofs of
\cite[Prop.~4.4]{BK01} and \cite[Lemma 3.2]{HKM06}. Indeed, we use the Feynman-Kac
formula for $U(t)$ and distinguish the contributions from those paths that leave
respectively do not leave the box $B_{t\log t}$ up to time $t$. The first
contribution can be estimated against the first term on the right of
\eqref{compact}, as is seen in \cite[Lemma 3.2]{HKM06}, together with the subsequent
text (see the display above (3.18) there). In order to see that the second
contribution can be estimated against the second term on the right hand side of
\eqref{compact}, we have to repeat parts of the proof of \cite[Prop.~4.4]{BK01}; we
shall replace the $R$ there by $2R\alpha(t)$.

As a first step, we estimate, with the help of a Fourier expansion, against the
principal eigenvalue. For $V\colon\Z^d\to\R$, let $\lambdaD_{t\log t}(V)$ be the
principal eigenvalue of $\DeltaD+V$ in the box $B_{t\log t}$ with zero boundary
condition. Then a Fourier expansion shows in a standard way that 
$$
\E_0\Big[{\rm e}^{\langle \ell_t,V\rangle}\1_{\{\supp(\ell_t)\subset B_{t\log
t}\}}\Big] 
\leq {\rm e}^{o(t/\alpha_t^{2})} {\rm e}^{t\lambdaD_{t\log t}(V)},
$$
where $o(t\alpha_t^{-2})$ does not depend on the potential $V$. Now
\cite[Prop.~4.4]{BK01} says that the eigenvalue in the box $B_{t\log t}$ may be
estimated from above against a small error plus the maximal eigenvalue in certain,
mutually overlapping boxes:
$$
\lambdaD_{t\log t}(V)\leq \max_{k\in B_{t\log
t}}\lambdaD_{4kR\alpha(t)+B_{3R\alpha(t)}}(V)+\frac {C}{R^2\alpha(t)^2},
$$
where $C>0$ does not depend on $V$ nor on $R$ nor on $t$. (Here $\lambdaD_{B}(V)$
denotes the eigenvalue of $\DeltaD+V$ in a bounded set $B\subset\Z^d$ with zero
boundary condition.) Hence, we may estimate
\begin{equation}\label{compact1}
\Big\langle{\rm e}^{t\lambdaD_{t\log
t}(\xi)}\1_{\widehat\Gamma_{t,\eps}}(\overline\xi_t)
\Big\rangle
\leq {\rm e}^{\frac{C t}{R^2 \alpha_t^2}}\sum_{k\in B_{t\log t}} \Big\langle{\rm
e}^{t\lambdaD_{4kR\alpha(t)+B_{3R\alpha(t)}}(\xi)}\1_{\widehat\Gamma_{t,\eps}}
(\overline\xi_t)\Big\rangle.
\end{equation}
Now we estimate $\1_{\widehat\Gamma_{t,\eps}}(\overline\xi_t)$. Observe that
\begin{equation}\label{Indikator}
\widehat\Gamma_{t,\eps}\subset \bigcap_{M\in(0,\infty)}\bigcap_{k\in B_{t\log
t}}\bigcap_{x\in Q_{2R}}\Big\{\psi\in\Fcal\colon\d_R\big(\psi(4kR+x+\cdot)\wedge
M,\widehat \psi(\cdot)\big)> \frac\eps2\Big\}.
\end{equation}
In order to show this, we show that the complement of the right side is contained in
the complement of the left side. Pick $k\in B_{t\log t}$ and $x\in Q_{2R}$ and
$\psi\in\Fcal$ such that 
$$
\frac\eps2\geq \d_R\big(\psi(4kR+x+\cdot)\wedge
M,\widehat\psi(\cdot)\big)=\int_{Q_R}\Big|{\rm e}^{\frac1 \rho(\psi(4kR+x+y)\wedge
M)}-{\rm e}^{\frac 1 \rho\widehat\psi(y)}\Big|\,\d y.
$$
Then, for $\widetilde x=4kR+x$, we have (recalling that $\phi(s)=\frac s{1+s}$ is
increasing in $s$), for any $M>0$,
$$
\begin{aligned}
\dist\big(\e^{\frac 1\rho(\psi(\widetilde x+\cdot)\wedge M)},\e^{\frac 1\rho\widehat
\psi(\cdot)}\big)
&\leq \sum_{r=1}^{R} 2^{-r}\phi\Big(\int_{Q_r}\Big|{\rm e}^{\frac1
\rho(\psi(\widetilde  x+y)\wedge M)}-{\rm e}^{\frac 1 \rho\widehat\psi(y)}\Big|\,\d
y\Big)+\sum_{r>R}2^{-r}\\
&\leq \phi\Big(\int_{Q_R}\Big|{\rm e}^{\frac1 \rho(\psi(\widetilde  x+y)\wedge
M)}-{\rm e}^{\frac 1 \rho\widehat\psi(y)}\Big|\,\d y\Big)+2^{-R}\leq
\frac\eps2+2^{-R}< \eps,
\end{aligned}
$$
by our assumption that $R>2-\log\frac \eps 2$. Hence, $\psi$ lies in
$\widehat\Gamma_{t,\eps}^{\rm c}$, which shows that \eqref{Indikator} holds. 

Now we use \eqref{Indikator} on the right hand side of \eqref{compact1} and obtain that
\begin{equation}\label{compact2}
\begin{aligned}
\Big\langle&{\rm e}^{t\lambdaD_{t\log
t}(\xi)}\1_{\widehat\Gamma_{t,\eps}}(\overline\xi_t)\Big\rangle\\
&\leq {\rm e}^{\frac{C t}{R^2 \alpha_t^2}}\sum_{k\in B_{t\log t}} \Big\langle{\rm
e}^{t \lambdaD_{4kR\alpha(t)+B_{3R\alpha(t)}}(\xi)}\prod_{\widetilde k\in B_{t\log
t}}\1_{\{\forall M>0\,\forall x\in Q_{2R}\colon\d_R(\overline\xi_t(4\widetilde
kR+x+\cdot)\wedge M,\widehat \psi(\cdot))>\frac\eps2\}}\Big\rangle\\
&\leq {\rm e}^{\frac{C t}{R^2 \alpha_t^2}}3^d (t\log t)^{d} \Big\langle{\rm e}^{t
\lambdaD_{3R\alpha(t)}(\xi)}\1_{\{\forall M>0\colon\overline\xi_t\wedge M\in
\Gamma_{R,\eps/2}\}}\Big\rangle,
\end{aligned}
\end{equation}
where we have estimated the product of indicators against the $k$-th factor, and we
have used the shift-invariance of the potential. Now enlarge $C$ in order to absorb
the term $3^d (t\log t)^{d}$.
\qed
\end{Proof}

\subsection{Truncating the potential.}\label{sec-potentialcutting}

In the next lemma, we replace the random potential by a truncated
version. In the proof of Lemma~\ref{Lemma_int} below it will turn out to be crucial
that the random potential under interest is bounded from above, hence
Lemma~\ref{lem-cutting} is a necessary preparation for that.

\begin{lemma}[Truncating the potential]\label{lem-cutting} Fix $R>0$ and $\eps>0$. Then
\begin{equation}\label{cuttingassertion}
\begin{aligned}
\limsup_{t\to\infty}&\frac{\alpha_t^2}t\log\Big\langle\E^{t,R}\Big[{\rm e}^{\langle
\ell_t, \xi_t\rangle}\Big]\1_{\{\forall M>0\colon\overline\xi_t\wedge M\in
\Gamma_{R,\eps}\}}\Big\rangle\\
&\leq
\limsup_{M\to\infty}\limsup_{t\to\infty}\frac{\alpha_t^2}t\log\Big\langle\E^{t,R}\Big[{\rm
e}^{\langle \ell_t, \xi_t\wedge
(M/\alpha_t^2)\rangle}\Big]\1_{\Gamma_{R,\eps}}(\overline\xi_t\wedge M)\Big\rangle.
\end{aligned}
\end{equation}
\end{lemma}
\begin{Proof}{Proof.} It is clear that we may estimate, for any $M>0$,
$$
\1_{\{\forall  \widetilde M>0\colon\overline\xi_t\wedge \widetilde M\in
\Gamma_{R,\eps}\}}\leq \1_{\Gamma_{R,\eps}}(\overline\xi_t\wedge M).
$$
Fix some small parameter $\eta>0$. In the expectation on the left hand side of
\eqref{cuttingassertion} we insert the sum of the indicators on the event $\{\langle
\ell_t,\xi_t-\xi_t\wedge (M/\alpha_t^2)\rangle\leq \eta t/\alpha_t^2\}$ and on the
opposite event. On the first event, we estimate $\langle \ell_t,\xi_t\rangle\leq
\eta t/\alpha_t^2+\langle \ell_t,\xi_t\wedge (M/\alpha_t^2)\rangle$ in the exponent.
The second indicator is estimated as follows:
$$
\1_{\{\langle \ell_t,\xi_t-\xi_t\wedge(M/\alpha_t^2)\rangle> \eta t/\alpha_t^2\}}
\leq \e^{-K\eta t/\alpha_t^2}\e^{K\langle \ell_t,\xi_t-\xi_t\wedge
(M/\alpha_t^2)\rangle},
$$
where $K\in(0,\infty)$ is some large auxiliary parameter. This gives that
\begin{equation}\label{cuttingesti}
\begin{aligned}
\Big\langle\E^{t,R}&\Big[{\rm e}^{\langle \ell_t, \xi_t\rangle}\Big]\1_{\{\forall 
\widetilde M>0\colon\overline\xi_t\wedge \widetilde M\in
\Gamma_{R,\eps}\}}\Big\rangle\\
&\leq \e^{\eta t/\alpha_t^2}\Big\langle\E^{t,R}\Big[{\rm e}^{\langle \ell_t,
\xi_t\wedge (M/\alpha_t^2)\rangle}\Big]\1_{\Gamma_{R,\eps}}(\overline\xi_t\wedge
M)\Big\rangle
+\e^{-K\eta t/\alpha_t^2}\Big\langle\E^{t,R}\Big[{\rm e}^{\langle \ell_t,
\xi_t\rangle}\e^{K\langle \ell_t,\xi_t-\xi_t\wedge
(M/\alpha_t^2)\rangle}\Big]\Big\rangle.
\end{aligned}
\end{equation}
The last expectation is estimated with the help of Schwarz' inequality:
\begin{equation}\label{CSineq}
\Big\langle\E^{t,R}\Big[{\rm e}^{\langle \ell_t, \xi_t\rangle}\e^{K\langle
\ell_t,\xi_t-\xi_t\wedge (M/\alpha_t^2)\rangle}\Big]\Big\rangle
\leq \Big\langle\E^{t,R}\Big[{\rm e}^{2\langle \ell_t,
\xi_t\rangle}\Big]\Big\rangle^{1/2}
\Big\langle\E^{t,R}\Big[\e^{2K\langle \ell_t,\xi_t-\xi_t\wedge
(M/\alpha_t^2)\rangle}\Big]\Big\rangle^{1/2}.
\end{equation}
We are going to show that, for any $K\in(0,\infty)$,
\begin{equation}\label{zerorate}
\limsup_{M\to\infty}\limsup_{t\to\infty}\frac{\alpha_t^2}t\log
\Bigl\langle\E^{t,R}\bigl[ \e^{K\langle \ell_t, \xi_t-(\xi_t\wedge
M/\alpha_t^2)\rangle}\bigr]\Bigr\rangle\leq 0,
\end{equation}
i.e., the exponential rate (on the scale $t/\alpha_t^2$) of the second term on the
right-hand side of \eqref{CSineq} vanishes as $M\to\infty$. In the course of the
proof, it will become obvious that the first term on the right-hand side of
\eqref{CSineq} has a bounded exponential rate. Hence, the assertion of the lemma
follows from considering the large-$t$ rate in \eqref{cuttingesti} and making
$M\to\infty$, $K\to\infty$ and $\eta\downarrow 0$.

We now prove \eqref{zerorate}. We first sum on all subsets $S$ of $B=B_{3R
{\alpha_t}}$ in which the potential $\xi_t$ is larger than $M/\alpha_t^2$ and
distinguish large and small such sets. This distinction is made with the help of a
small auxiliary parameter $\tau\in(0,\infty)$:
\begin{equation}\label{A_beide}
\begin{aligned}
\Bigl\langle\E^{t,R}\bigl[ \e^{K\langle \ell_t, \xi_t-(\xi_t\wedge
M/\alpha_t^2)\rangle}\bigr]\Bigr\rangle 
&\le  \Bigl\langle\E^{t,R} \Bigl[\exp\Big\{K\sum_{z \in
B}\ell_t(z)\cdot\xi_t(z)\cdot \1_{\{\xi_t(z)>M/\alpha_t^2\}}\Big\}
\Bigr]\Bigr\rangle \\
&\le \sum_{S \subset B \colon |S|\ge \tau \alpha_t^d}\Bigl\langle\E^{t,R}
\Bigl[\e^{K\sum_{z \in S}\ell_t(z)\xi_t(z)} \Bigr]\1_{\{S=\{z\in B\colon
\xi_t(z)>M/\alpha_t^2\}\}}\Bigr\rangle \\
&\qquad + \sum_{S \subset B \colon |S|< \tau \alpha_t^d}\Bigl\langle\E^{t,R}
\Bigl[\e^{K\sum_{z \in S}\ell_t(z)\cdot\xi_t(z)} \Bigr]\Bigr\rangle.
\end{aligned}
\end{equation}
In the following, we show that the exponential rate of the first term tends to
$-\infty$ as $M\to\infty$ for any $\tau>0$, and the rate of the second vanishes as
$\tau\downarrow 0$.

We first consider a summand of the first sum, where $|S|\ge \tau \alpha_t^d$. By
using Schwarz' inequality we obtain
\begin{equation}\label{xi_ab_A}
\begin{aligned}
\Bigl\langle&\E^{t,R} \Bigl[\e^{K\sum_{z \in S}\ell_t(z)\cdot\xi_t(z)}
\Bigr]\1_{\{S=\{z\in B\colon \xi_t(z)>M/\alpha_t^2\}\}}\Bigr\rangle \\
&\le  \Bigl\langle \E^{t,R} \Bigl[ \e^{2K\sum_{z \in S}\ell_t(z)\cdot\xi_t(z)} 
\Bigr] \Bigr\rangle^\frac{1}{2}\cdot \Prob\Bigl(S=\bigr\{z\in
B\colon\xi_t(z)>M/\alpha_t^2\bigl\}\Bigr)^\frac{1}{2} 
\end{aligned}
\end{equation}
We first estimate the probability. We use the definition of $\xi_t(z)$ in
\eqref{xitdef}, the independence of the $\xi(z)$ for different $z$, the Markov
inequality and the definition of $H(\cdot)$ in (\ref{Def_H}), to obtain
\begin{equation}\label{Prob}
\begin{aligned}
\Prob\Bigl(S=\bigr\{z\colon\xi_t(z)>M/\alpha_t^2\bigl\}\Bigr)^\frac{1}{2}
&\le \Prob\Bigl(\xi(0)
>\frac{H(t/\alpha_t^{d})}{t/\alpha_t^{d}}+\frac{M}{\alpha_t^2} \Bigr)^\frac{|S|}{2}
\\
&=\Prob\Bigl(\e^{\xi(0)t/\alpha_t^d} > \e^{H(t/\alpha_t^{d})}
\e^{Mt/\alpha_t^{d+2}}\Bigr)^\frac{|S|}{2}\\
&\le\Bigl( \e^{-M t/\alpha_t^{d+2}} \e^{-H(t /\alpha_t^{d})} 
\Bigl\langle \e^{\xi(0)t/\alpha_t^d}\Bigr\rangle \Bigr)^\frac{|S|}{2}\\
&\le \e^{-t\alpha_t^{-2}M  \tau/2}.
\end{aligned}
\end{equation}
The exponential rate of this tends to $-\infty$ as $M\to\infty$. Hence, it suffices
to show that the exponential rate of  the expectation on the right hand side of
(\ref{xi_ab_A}) is finite on the scale $t\alpha_t^{-2}$. This is surprisingly
difficult and cannot be handled with rough arguments. We do this by extending some
results of \cite[Sect. 3.2]{HKM06}, the only additional issue being that the set $B$
is replaced by some subset $S$ of $B$. This is some technical issue since $\xi_t$
can assume also negative values, such that we have to repeat some of the steps from
\cite{HKM06}. We use the definition of $\xi_t$ in \eqref{xitdef}, apply Fubini's
theorem, execute the expectation with respect to $\xi$, recall the definition of
$H(\cdot)$ in (\ref{Def_H}) and  use the abbreviation
$h_t(z)=H(2K\ell_t(z))-2K\ell_t(z)\frac{H(t/\alpha_t^{d})}{t/\alpha_t^{d}}$. This
gives
\begin{equation}\label{Sum_1}
\begin{aligned}
\Bigl\langle \E^{t,R} \Big[ \e^{2K\sum_{z \in S}\ell_t(z)\cdot\xi_t(z)}  \Big]
\Bigr\rangle
&=\E^{t,R} \Big[ \e^{\sum_{z \in S} h_t(z)}  \Big].
\end{aligned}
\end{equation}
We split the sum on $z\in S$ into the subsums where $\ell_t(z)\leq t/\alpha_t^d$ and
the remainder. For $\ell_t(z)\leq t/\alpha_t^d$ we may apply the asymptotics for $H$
from Lemma~\ref{alpha}. This gives, as $t\to\infty$, also using the relation between
$\alpha_t$ and $\kappa(t)$ in \eqref{Def_alpha},
\begin{equation}\label{Sum_2}
\begin{aligned}
\sum_{z \in S\colon \ell_t(z)\leq t/\alpha_t^d} h_t(z)
&\leq (\rho+o(1)) \kappa(t/\alpha_t^d) \sum_{z \in S\colon \ell_t(z)\leq
t/\alpha_t^d
}2K\ell_t(z)\tfrac{\alpha_t^d}{t}\log\Big(2K\ell_t(z)\tfrac{\alpha_t^d}{t}\Big)\\
&\leq 2\rho 2K\log(2K)\frac t{\alpha_t^{d+2}}|S|\leq C \frac {|S|}{|B|}\frac
t{\alpha_t^{2}},
\end{aligned}
\end{equation}
where $C\in(0,\infty)$ depends on $\rho$ and $K$ (and $R$) only.

Now we handle the subsum on $z\in S$ satisfying $\ell_t(z)>t/\alpha_t^d$. For this
purpose, we need the following estimate for differences of $H$-terms, which follows from
\cite[Theorem~3.8.6.(a)]{BGT87}. For any $\delta\in(0, \frac 12]$, there are
$A,t_0\in(1,\infty)$ such that
\begin{equation}\label{Hupperesti}
\frac{H(ty)-yH(t)}{\kappa(t)}\leq A y^{1+\delta},\qquad y\in [1,\infty),
t\in[t_0,\infty). 
\end{equation}
We pick a small $\delta>0$ and apply \eqref{Hupperesti} with $\delta$ replaced by
$\delta^2/3$, for $y=2K\ell_t(z)\alpha_t^d/t$ and with $t/\alpha_t^{d}$ instead of
$t$, to get, for $z$ satisfying $\ell_t(z)>t/\alpha_t^d$,
\begin{equation}\label{Absch_H}
\begin{aligned}
h_t(z)=H(2K \ell_t(z))-2K\ell_t(z) \frac{\alpha_t^d}{t} H( t/\alpha_t^d)
&\leq A \kappa(t/\alpha_t^{d})
\Big(2K\ell_t(z)\tfrac{\alpha_t^d}{t}\Big)^{1+\delta^2/3}\\
&= \widetilde C \frac{t}{\alpha_t^{d+2}}
\Big(\frac{\alpha_t^d}{t}\ell_t(z)\Big)^{1+\delta^2/3},
\end{aligned}
\end{equation}
where we have used the definition of $\alpha_t$  in (\ref{Def_alpha}), and
$\widetilde C\in(0,\infty)$ depends on $A$ and $K$ only. Using this and
\eqref{Sum_2} in \eqref{Sum_1}, we can estimate
\begin{equation}\label{Absch_H1}
\E^{t,R} \Big[ \e^{\sum_{z \in S} h_t(z)}  \Big]
\leq \e^{C \frac{|S|}{|B|}t/\alpha_t^2}\,\E^{t,R} \Big[ \exp\Big\{\widetilde C
\frac{t}{\alpha_t^{2}} \alpha_t^{d\delta^2/3}\sum_{z \in S\colon
\ell_t(z)>t/\alpha_t^d}\big(\sfrac 1t\ell_t(z)\big)^{1+\delta^2/3}\Big\}\Big].
\end{equation}
Now in the same way as in \cite[Sect. 3.2]{HKM06}, we see that, for any probability
measure $\mu$ on $B$, for any $0<a\leq b<1$ and $0<c$,
\begin{equation}\label{Trick}
\sum_{z\in B\colon \mu(z)>\alpha_t^{-d}}\mu(z)^{1+a}\leq
\alpha_t^{d[(1-b)(b-a)+c(1+a-b)]}\Big(\sum_{z\in B}\mu(z)^{1+b+c}\Big)^{1+a-b}.
\end{equation}
This is proved as follows, using Jensen's inequality (we write $\sum_z$ instead of
$\sum_{z\in B\colon \mu(z)>\alpha_t^{-d}}$):
$$
\begin{aligned}
\sum_{z}\mu(z)^{1+a}
&=\Big(\sum_{z}\mu(z)^b\Big)\Big(\sum_{z}\frac{\mu(z)^b}{\sum_{z}\mu(z)^b}
\mu(z)^{1+a-b}\Big)\\
&\leq \Big(\sum_{z}\mu(z)^b\Big)
\Big(\sum_{z}\frac{\mu(z)^{1+b}}{\sum_{z}\mu(z)^b}\Big)^{1+a-b}\\
&\leq
\Big(\sum_{z}\mu(z)^b\big(\mu(z)\alpha_t^d\big)^{1-b}\Big)^{b-a}\Big(\sum_{z}\mu(z)^{1+b}\big(\mu(z)\alpha_t^d\big)^c\Big)^{1+a-b}\\
&\leq \alpha_t^{d[(1-b)(b-a)+c(1+a-b)]}\Big(\sum_{z\in B}\mu(z)^{1+b+c}\Big)^{1+a-b}.
\end{aligned}
$$

Applying \eqref{Trick} for $\mu=\frac 1t\ell_t$, $a=\delta^2/3$, $b=\delta^2/3+\frac
\delta{1+\delta}$ and $c=\delta-\delta^2/3-\frac \delta{1+\delta}$, we obtain
\begin{equation}
\frac{t}{\alpha_t^{2}}\alpha_t^{d\delta^2/3}\sum_{z\in B\colon
\ell_t(z)>t/\alpha_t^d}\big(\sfrac 1t\ell_t(z)\big)^{1+\delta^2/3}
\leq  \alpha_t^{-\frac 1{1+\delta}(d+(2-d)(1+\delta))}\|\ell_t\|_{1+\delta},
\end{equation}
where $\|\cdot\|_{1+\delta}$ denotes the $(1+\delta)$-norm on
$\ell^{1+\delta}(\Z^d)$. Now, picking $\delta>0$ so small that $\delta(d-2)<2$,
\cite[Prop.~2.1]{HKM06} states that the large-$t$ exponential rate of the right hand
side of \eqref{Absch_H1} on the scale $t/\alpha_t^2$ vanishes as $\widetilde
C\downarrow 0$. However, the proof shows that this rate is finite for any
$\widetilde C\in(0,\infty)$. Using this fact in \eqref{Sum_1}, and substituting this
in \eqref{xi_ab_A} we see because of \eqref{Prob} that the exponential rate of the
first sum on the right hand side of \eqref{A_beide}  on the scale $t/\alpha_t^2$
tends to $-\infty$ as $M\to\infty$, for any $\tau>0$.

Now we address the second sum on the right hand side of (\ref{A_beide}). We show
that its large-$t$ exponential rate on the scale $t/\alpha_t^2$ vanishes as
$\tau\downarrow 0$. We consider $S\subset B=B_{3R\alpha_t^d}$ with
$|S|<\tau\alpha_t^d$. We start from \eqref{Absch_H1}, which is valid for any
$S\subset B$. We use H\"older's inequality with new parameters $\frac 1p+\frac 1q=1$
for the last sum to obtain
\begin{equation}\label{Trick2}
\sum_{z \in S\colon \ell_t(z)>t/\alpha_t^d}\big(\sfrac 1t\ell_t(z)\big)^{1+\delta^2/3}
\leq |S|^{\frac 1q}\Big(\sum_{z \in B\colon \ell_t(z)>t/\alpha_t^d}\big(\sfrac
1t\ell_t(z)\big)^{p(1+\delta^2/3)}\Big)^{\frac 1p}.
\end{equation}
Now we apply \eqref{Trick} for $\mu=\frac 1t\ell_t$, $a=p+p\delta^2/3-1$, some $b\in
(a,1)$ and $c=\frac{p-(1+b)(1+a-b)}{1+a-b}$, where we assume that $\delta>0$ is
small enough and $p>1$ close enough to one such that all the assumptions $0<a\leq
b<1$ and $c>0$ are satisfied. This gives, using \eqref{Trick2},
$$
\begin{aligned}
\frac{t}{\alpha_t^{2}}&\alpha_t^{d\delta^2/3}\sum_{z\in S\colon
\ell_t(z)>t/\alpha_t^d}\big(\sfrac 1t\ell_t(z)\big)^{1+\delta^2/3}\\
&\leq \frac{t}{\alpha_t^{2}}\alpha_t^{d\delta^2/3}\tau^{1/q}\alpha_t^{d/q}\,
\alpha_t^{d[(1-b)(b-a)+p-(1+b)(1+a-b)]/p}\|\sfrac 1t\ell_t\|_{1+b+c}\\
&=\tau^{1/q}\alpha_t^{-\frac
1{1+\widetilde\delta}(d+(2-d)(1+\widetilde\delta))}\|\ell_t\|_{1+\widetilde\delta},
\end{aligned}
$$
where $\widetilde\delta=b+c$. Picking $\delta>0$ small enough and $p>1$ close
enough to one, we also have that $\widetilde\delta(d-2)<2$, and we may again
apply \cite[Prop.~2.1]{HKM06} and see that the exponential rate of the second sum on
the right hand side of (\ref{A_beide}) vanishes as $\tau\downarrow0$. This ends the
proof.
\qed\end{Proof}

\subsection{Main part of the proof of
Proposition~\ref{Hauptprop}.}\label{sec-PropProof}

Using Lemmas~\ref{Compact} and \ref{lem-cutting}, it is clear that
Proposition~\ref{Hauptprop} now follows from the following assertion.

\begin{prop}\label{ZweiteProp} For any $\eps>0$,
\begin{equation}\label{mainstep}
\limsup_{R\to\infty}\limsup_{M\to\infty}\limsup_{t\to\infty}\frac{\alpha_t^2}t\log
\Big\langle\E^{t,R}\Big[{\rm e}^{\langle \ell_t, \xi_t\wedge
(M/\alpha_t^2)\rangle}\Big]\1_{\Gamma_{R,\eps}}(\overline\xi_t\wedge M)\Big\rangle
<-\chi.
\end{equation}
\end{prop}

Let us now prove Proposition~\ref{ZweiteProp}. For any potential $V\colon B_R\to\R$,
we denote by $\lambdaD_R(V)$ the principal eigenvalue of $\DeltaD+V$ in the box
$B_{R}$ with zero boundary condition. Introduce a rescaled version of this
eigenvalue by putting, for $\psi\in\Fcal(Q_R)$,
\begin{equation}\label{rescaledeigenv}
\lambda_{R}^{\ssup t}(\psi)=\alpha_t^2 \lambda^{\rm
d}_{R\alpha(t)}\big(\sfrac{1}{\alpha_t^2}\psi^{\rm
d}\big),\qquad\mbox{where}\quad\psi^{\rm d}(z)=\alpha_t^d
\int_{z/\alpha_t+[0,\alpha_t^{-1})^d} \psi(y)\,\d y\quad \mbox{for }z \in \Z^d.
\end{equation}
Observe from \eqref{xitdef} and \eqref{xitbardef} that
$\sfrac{1}{\alpha_t^2}\overline\xi_t^{\rm d}=\xi_t$. Recall the definition of
$\E^{t,R}$ from Lemma~\ref{Compact}. Hence, using a Fourier expansion, one has, for
any $R,M>0$ , as $t\to\infty$,
\begin{equation}\label{eigenvalue}
\begin{aligned}
\E^{t,R}\Big[{\rm e}^{\langle \ell_t, \xi_t\wedge (M/\alpha_t^2)\rangle}\Big]
&= \e^{o(t/\alpha_t^2)} \exp\Big\{t \lambda^{\rm d}_{3R\alpha_t}\big(\xi_t\wedge
(M/\alpha_t^2)\big)\Big\}\\
&=\e^{o(t/\alpha_t^2)}\exp\Big\{\frac t{\alpha_t^2} \lambda_{3R}^{\ssup
t}(\overline\xi_t\wedge M)\Big\}.
\end{aligned}
\end{equation}
Now we multiply \eqref{eigenvalue} with $\1_{\Gamma_{R,\eps}}(\overline\xi_t\wedge
M)$ and take expectation with respect to  $\xi$. We subtract and add the term
$t\alpha_t^{-2}\rho\log(\frac \e\rho\Lcal_{3R}(\overline \xi_t\wedge M))$ in the
exponent. The next step is to pick some
small parameter $\beta\in(0,\infty)$ and to distinguish the events $\{\overline
\xi_t\wedge M\in D_{\beta,R}\}$ and its complement, where
\begin{equation}\label{Dbetadef}
D_{\beta,R}=\big\{\psi\in\Fcal(Q_{3R})\colon|\Lcal_{3R}(\psi)-\rho|\leq \beta\big\}.
\end{equation}
On the event $\{\overline\xi_t\wedge M\in\Gamma_{R,\eps}\}$, we estimate the first
two terms in the exponent differently on $\{\overline \xi_t\wedge M\in
D_{\beta,R}\}$ and on the complement:
\begin{equation}
\lambda_{3R}^{\ssup t}(\overline\xi_t\wedge M)-\rho\log\Big(\frac
\e\rho\Lcal_{3R}(\overline \xi_t\wedge M)\Big)
\leq\begin{cases}-\chi_R(\beta,\eps,t)&\mbox{on }\{\overline \xi_t\wedge M\in
D_{\beta,R}\},\\
-\chi_{R}(t)&\mbox{on }\{\overline \xi_t\wedge M\notin D_{\beta,R}\},
\end{cases}
\end{equation}
where the variational formulas are defined by
\begin{equation}\label{chiapproxdef}
\begin{aligned}
-\chi_R(\beta,\eps,t)&=\sup\Big\{\lambda_{3R}^{\ssup t}(\psi)-\rho\log \Big(\frac
\e\rho\Lcal_{3R}(\psi)\Big)\colon \psi\in \Gamma_{R,\eps}\cap D_{\beta,R}\Big\},
\end{aligned}
\end{equation}
and we put
$-\chi_{R}(t)=-\chi_{R}(\infty,0,t)=\sup_{\psi\in\Fcal(Q_{3R})}[\lambda_{3R}^{\ssup
t}(\psi)-\rho\log (\frac \e\rho\Lcal_{3R}(\psi))] $.

Making the above explicit and abbreviating
\begin{equation}\label{Ftdef}
F_{t,R}(\psi)=\exp\Big\{\frac t{\alpha_t^2}\rho\log\Big(\frac
\e\rho\Lcal_{3R}(\psi)\Big)\Big\},
\end{equation}
we obtain
\begin{equation}\label{kernelcalc}
\begin{aligned}
\Big\langle&\E^{t,R}\Big[{\rm e}^{\langle \ell_t, \xi_t\wedge
(M/\alpha_t^2)\rangle}\Big]\1_{\Gamma_{R,\eps}}(\overline\xi_t\wedge M)\Big\rangle\\
&\leq \e^{o(t/\alpha_t^2)}\Big\langle\exp\Big\{\frac t{\alpha_t^2}\Big[
\lambda_{3R}^{\ssup t}(\overline\xi\wedge M)-\rho\log\Big(\frac
\e\rho\Lcal_{3R}(\overline \xi_t\wedge M)\Big)\Big]\Big\} F_{t,R}(\overline
\xi_t\wedge M) \1_{\Gamma_{R,\eps}}(\overline\xi_t\wedge M)\Big\rangle\\
&\leq \e^{o(t/\alpha_t^2)}\Big\langle\Big[ \1_{D_{\beta,R}}(\overline \xi_t\wedge
M)\e^{-t\alpha_t^{-2}\chi_R(\beta,\eps,t)}+\1_{D_{\beta,R}^{\rm c}}(\overline
\xi_t\wedge M) \e^{-t\alpha_t^{-2}\chi_R(t)}\Big]\\
&\qquad\qquad \qquad\qquad \times F_{t,R}(\overline \xi_t\wedge M)
\1_{\Gamma_{R,\eps}}(\overline\xi_t\wedge M)\Big\rangle \\
&\leq \e^{o(t/\alpha_t^2)}\Big[\e^{-t\alpha_t^{-2}\chi_R(\beta,\eps,t)}
\Big\langle F_{t,R}(\overline \xi_t\wedge M)\Big\rangle
+\e^{-t\alpha_t^{-2}\chi_R(t)}\Big\langle F_{t,R}(\overline \xi_t\wedge M)
\1_{D_{\beta,R}^{\rm c}}(\overline \xi_t\wedge M)\Big\rangle\Big].
\end{aligned}
\end{equation}

Now we need the following asymptotics for the approximative variational formulas:

\begin{lemma}\label{lem-apprvarform}
\begin{itemize}
\item[(i)]
$$
\liminf_{R\to\infty}\liminf_{t\to\infty}\chi_R(t)\geq \chi.
$$
\item[(ii)] For any $\eps>0$, and any $\beta>0$ small enough,
$$
\liminf_{R\to\infty}\liminf_{t\to\infty}\chi_R(\beta,\eps,t)> \chi.
$$
\end{itemize}
\end{lemma}

The proof of Lemma~\ref{lem-apprvarform} is deferred to the end of
Section~\ref{sec-VarForm}.

The large-$t$ exponential rate of the expectation of $F_{t,R}(\overline \xi_t\wedge
M)$ is nonpositive for any $M>0$, as is seen from an application of part (i) of the
following lemma for $K=\rho$.

\begin{lemma}\label{Lemma_int}
\begin{itemize}
\item[(i)] Fix $R>0$ and $M>0$. Then, for any $K>0$,
\begin{equation}\label{upperboundlogL}
\limsup_{t \to \infty}\frac{\alpha_t^2}{t} \log \Big\langle F_{t,R}(\overline
\xi_t\wedge M)^{K/\rho}\Big\rangle \le K\log \frac {K}\rho.
\end{equation}
\item[(ii)] For any $\beta>0$,
\begin{equation}
\limsup_{R \to \infty}\limsup_{M \to \infty}\limsup_{t \to
\infty}\frac{\alpha_t^2}{t} \log \Big\langle F_{t,R}(\overline \xi_t\wedge M)
\1_{D_{\beta,R}^{\rm c}}(\overline \xi_t\wedge M)\Big\rangle <0.
\end{equation}
\end{itemize}
\end{lemma}

It is elementary to see that an application of Lemmas~\ref{lem-apprvarform} and
\ref{Lemma_int} to the terms on the right hand side of \eqref{kernelcalc} implies
that \eqref{mainstep} holds for any $\eps>0$. This ends the proof of
Proposition~\ref{ZweiteProp}.

It remains to prove Lemma~\ref{Lemma_int}. Let us mention that our proof crucially
depends on the appearance of the cut-off potential $\overline\xi_t\wedge M$ instead
of $\overline\xi_t$, even though the cut-off parameter $M$ does not appear in the
asymptotics. This is the place where Lemma~\ref{lem-cutting} turns out to be
important.

\begin{Proof}{Proof of Lemma~\ref{Lemma_int}.}
First we prove (i). Recall the definition of $F_{t,R}(\psi)$ from \eqref{Ftdef}. We
recall that $\overline \xi_t(x)=\alpha_t^2\xi_t(\lfloor x\alpha_t\rfloor)$ and
rewrite 
\begin{equation}\label{lem-int1}
\Big\langle F_{t,R}(\overline \xi_t\wedge M)^{K/\rho} \Big\rangle
=\Big\langle \Bigl( \alpha_t^{-d} \sum\limits_{z \in B}
\exp\Big\{\frac1\rho\big(\alpha_t^2 \xi_t(z)\wedge M\big)\Big\} \Bigr)^{D_t}
\Big\rangle,
\end{equation}
where we abbreviated $D_t=Kt/\alpha_t^2$ and $B=B_{3R\alpha_t}$. (For the ease of
notation, we assume that $D_t$ and $3R\alpha_t$ are integers.)

Now we calculate the right hand side with the help of elementary combinatorics. We
denote by $\mathcal{M}_1^{\ssup{D_t}}(B)\!=\!\bigl\{\mu \in
(\tfrac{1}{D_t}\mathbb{N}_0)^B\colon\sum_{b \in B}\mu(b)=1\bigr\}$ the set of
probability vectors $\mu$ on $B$ such that $D_t\mu$ has throughout integer
coefficients. Then we
have
\begin{equation}\label{lem-int2} 
\begin{aligned}
\mbox{l.h.s.~of }\eqref{lem-int1}
&=\alpha_t^{-d  D_t}\sum\limits_{z_1,z_2,\ldots,z_{D_t} \in B} \biggl\langle
\prod\limits_{b \in B} \exp\Bigl\{ \Bigl(\sfrac{\alpha_t^2}{\rho}\xi_t(b) \wedge
\sfrac{M}{\rho}\Bigr) \#\{i\colon z_i=b\}\Bigr\} \biggr\rangle\\
&=\alpha_t^{-d  D_t}\sum_{\mu \in \mathcal{M}_1^{\ssup{D_t}}(B)}\#\Bigl\{z \in
B^{D_t}\colon \# \{i\colon z_i=b\}=D_t \mu(b), \forall \,b\Bigr\}\\
&\qquad\qquad\qquad\qquad\times
\prod_{b \in B}\biggl\langle  \e^{\frac K\rho t \mu(b) \bigl(\xi_t(0)\wedge
\frac{M}{\alpha_t^2}\bigr)} \biggr\rangle \\
&=\alpha_t^{-d  D_t}\sum_{\mu \in
\mathcal{M}_1^{\ssup{D_t}}(B)}\frac{D_t!}{\prod\limits_{b \in B} (D_t
\mu(b))!}\;\prod\limits_{b \in B}\biggl\langle  \e^{\frac K\rho t \mu(b)\bigl(
\xi_t(0)\wedge \sfrac{M}{\alpha_t^2}\bigr)}\biggr\rangle.
\end{aligned}
\end{equation}
Use Stirling's formula and recall that $D_t=Kt/\alpha_t^2$ to deduce that, uniformly
in $\mu \in \mathcal{M}_1^{\ssup{D_t}}(B)$,
\begin{equation}\label{Stirling}
\alpha_t^{-d  D_t}\frac{D_t!}{\prod_{b \in B} (D_t
\mu(b))!}=\e^{o(t/\alpha_t^2)}\exp\Big\{-K\frac t{\alpha_t^2}\sum_{b\in
B}\mu(b)\log(\mu(b)\alpha_t^d)\Big\}.
\end{equation}

Now we analyse the  last product on the right hand side of (\ref{lem-int2}). 
We use the formula $\E[X]=\int_0^\infty \P(X>s)\,\d s$ for nonnegative random
variables $X$, introduce an auxiliary variable $N>0$ (which will be chosen later)
and apply the Markov inequality with the map $s\mapsto s^N$. Then we make the change
of measure via $s=\exp\{\frac{t}{\alpha_t^{d+2}} r\}$, which implies $\tfrac{\d
s}{\d r}=\tfrac{t}{\alpha_t^{d+2}}\exp\{\frac{t}{\alpha_t^{d+2}} r\}$. We use the
abbreviation $a=\frac K\rho\mu(b)\alpha_t^d$. Hence, we have, for any $Q\geq 0$,
\begin{equation}\label{lem-int3}
\begin{aligned}
\biggl\langle \e^{\frac K\rho t \mu(b)\bigl(\xi_t(0)\wedge
\sfrac{M}{\alpha_t^2}\bigr)}\biggr\rangle
&\le \e^{tQ/\alpha_t^{d+2}}+\int^{\exp\{tMa/\alpha_t^{d+2}\}}_{\exp\{tQ/\alpha_t^{d+2}\}}
\text{Prob}\Bigl[\e^{\frac{t}{\alpha_t^d}\xi_t(0) a} >s \Bigr] \, \d s\\
&\leq  \e^{tQ/\alpha_t^{d+2}}+\int_{\exp\{tQ/\alpha_t^{d+2}\}}^{\exp\{tMa/\alpha_t^{d+2}\}} s^{-N}\Big\langle
\exp\Big\{ \frac{t}{\alpha_t^{d}} \xi_t(0) a N \Big\}\Big\rangle  \, \d s \\
&= \e^{tQ/\alpha_t^{d+2}}+\frac{t}{\alpha_t^{d+2}}\int_Q^{Ma } \exp\Big\{r(1-N)
\frac{t}{\alpha_t^{d+2}}\Big\} \Big\langle \exp\Big\{\frac{t}{\alpha_t^{d}} \xi_t(0)
a N \Big\}\Big\rangle  \, \d r .
\end{aligned}
\end{equation}
Now we use the definition $\xi_t(z)=\xi(z)-H(t/\alpha_t^d)\alpha_t^d/t$ of the
shifted potential (see \eqref{xitdef}) and recall that  $H(s)=\log\langle
\e^{s\xi(0)}\rangle$ and $t/\alpha_t^{d+2}=\kappa(t/\alpha_t^d)$, to proceed with 
\begin{equation*}
\begin{aligned}
\mbox{l.h.s.~of }\eqref{lem-int3}
&\le \e^{tQ/\alpha_t^{d+2}}+\frac{t}{\alpha_t^{d+2}}\int_Q^{Ma }
\exp\Big\{-(N-1)r\frac{t}{\alpha_t^{d+2}}\Big\}\\
&\qquad\qquad\qquad\qquad\times\exp\Big\{-aN
H\Big(\frac{t}{\alpha_t^d}\Big)\Big\} \Big\langle \exp\Big\{\frac{t}{\alpha_t^{d}}
\xi(0) a N \Big\}\Big\rangle  \,\d r \\
&= \e^{tQ/\alpha_t^{d+2}}+\frac{t}{\alpha_t^{d+2}} \int_Q^{Ma } \exp\Big\{-\frac{t}{\alpha_t^{d+2}}
\Big((N-1)r+\frac{aN H(t/\alpha_t^d)-H(aN t/\alpha_t^d)}{\kappa(t/\alpha_t^d)}
\Big)\Big\} \, \d r.
\end{aligned}
\end{equation*}
Now we have to distinguish the case of bounded $aN$, where we can use precise
asymptotics in (\ref{Ann_H}) for the last quotient, and the case of arbitrarily
large $aN$, where we can only bound the last quotient. Introduce a new parameter
$L>0$, which will later be chosen large enough. First we handle those $a$ satisfying
$a\leq L$, and we now pick $N=\e^{r/a \rho}/a\e$. Note that $aN$ lies then in the interval $
[\tfrac{1}{e},\tfrac{e^{M/\rho}}{e}]$.
Hence, we may use the asymptotics in (\ref{Ann_H}). This gives, picking $Q=0$,
\begin{equation*}
\begin{aligned}
\mbox{l.h.s.~of }\eqref{lem-int3}
&\le 1+ \e^{o(t/\alpha_t^{d+2})} \int_0^{Ma } \exp\Big\{-\frac{t}{\alpha_t^{d+2}}
\Bigl((N-1)r-\rho a N \log(aN) \Bigr)\Bigr\} \, \d r \\
&= 1+ \e^{o(t/\alpha_t^{d+2})} \int_0^{Ma } \exp\Bigl\{-\frac{t}{\alpha_t^{d+2}} 
\Bigl(-r+\sfrac{\rho}{\e}\e^\frac{r}{a \rho} \Bigr)\Bigr\} \, \d r.
\end{aligned}
\end{equation*}
The term $-r+\sfrac{\rho}{\e}\e^\frac{r}{a \rho}$ is  minimal for $r=\rho a \log(a 
\e)$ with value $-\rho a \log(a)$. Hence, we may estimate for $a\leq L$ as follows.
\begin{equation}\label{Fall_1a}
\begin{aligned}
\mbox{l.h.s.~of }\eqref{lem-int3}
&\le 2+ \e^{o(t/\alpha_t^{d+2})} M a \exp\Bigl\{\frac{t}{\alpha_t^{d+2}} \rho a 
\log(a) \Bigr\}
& \le \e^{o(t/\alpha_t^{d+2})}\exp\Bigl\{\frac{t}{\alpha_t^{d+2}} \rho a  \log(a)
\Bigr\}.
\end{aligned}
\end{equation}

Now we turn to $a$ satisfying $a>L$. This time we pick
$N=r^{1/\delta}(A(1+\delta))^{-1/\delta}a^{-(1+\delta)/\delta}$, where we have
picked some small $\delta >0$. With $A$ as in \eqref{Hupperesti} we have, for every
$t$ large enough,
\begin{equation*}
\frac{aN H(t/\alpha_t^d)-H(aNt/\alpha_t^d)}{\kappa(t/\alpha_t^d)}\ge -A\cdot
(aN)^{1+\delta}, \qquad  aN \ge 1.
\end{equation*}
This time we pick $Q=A(1+\delta)L$ and note that $aN\geq 1$ on the integration interval $[A(1+\delta)L, Ma]$. Hence, for $a>L$, we may estimate
$$
\mbox{l.h.s.~of }\eqref{lem-int3}
\le \exp\Bigl\{\frac{A(1+\delta)L\cdot t}{\alpha_t^{d+2}}\Bigl\}+
\e^{o(t/\alpha_t^{d+2})} \int_{A(1+\delta)L}^{Ma }
\exp\Bigl\{\frac{t}{\alpha_t^{d+2}}\Bigl(r-Nr+A(Na)^{1+\delta} \Bigr)\Bigr\} \, \d r.
$$
Note that we may estimate $-Nr+A(Na)^{1+\delta}\leq 0$ in the exponent. Furthermore we extend the integration area to the interval $[0,Ma]$. Hence,
\begin{equation}\label{Fall_2a}
\begin{aligned}
\mbox{l.h.s.~of }\eqref{lem-int3}
&\le \exp\Bigl\{\frac{A(1+\delta)L\cdot t}{\alpha_t^{d+2}}\Bigl\}+\e^{o(t/\alpha_t^{d+2})} \int_0^{Ma }\exp\Big\{\frac{t}{\alpha_t^{d+2}}  r \Big\}
\,\d r\\
&\le \e^{o(t/\alpha_t^{d+2})}\exp\Big\{\frac{t}{\alpha_t^{d+2}} Ma\Big\},
\end{aligned}
\end{equation}
where the last step is valid for $M>A(1+\delta)$.

Now we go back to \eqref{lem-int2} and substitute \eqref{Stirling}, recall that
$a=\frac K\rho \mu(b)\alpha_t^d$ and substitute \eqref{Fall_1a} for $a\leq L$ and
\eqref{Fall_2a} for $a>L$. We now write $L\rho/K$ instead of $L$ and obtain
\begin{equation}\label{Produkt}  
\begin{aligned}
\mbox{l.h.s.~of }\eqref{lem-int1}
&\le \e^{o(t/\alpha_t^2)}  \sum_{\mu \in \mathcal{M}_1^{\ssup{D_t}}(B)}
\Big(\prod_{b \in B}\exp\Bigl\{-\frac{t}{\alpha_t^2}K
\mu(b)\log\bigl(\mu(b)\alpha_t^d\bigr) \Bigr\}\Big) \\
&\qquad\times\Big(\prod_{b \in B \colon\mu(b)\alpha_t^d\le L}
\exp\Bigl\{\frac{t}{\alpha_t^{d+2}} K \mu(b)\alpha_t^d  \log\Big(\frac K\rho
\mu(b)\alpha_t^d\Big)\Bigr\}\Big)\\
&\qquad\times
\Big(\prod_{b \in B \colon \mu(b)\alpha_t^d > L}\exp\Big\{\frac{t}{\alpha_t^{d+2}}
M\frac K\rho\mu(b)\alpha_t^d\Big\}\Big)\\
&= \e^{o(t/\alpha_t^2)}   \sum_{\mu \in \mathcal{M}_1^{\ssup{D_t}}(B)}
\exp\Bigl\{\frac{t}{\alpha_t^2}K\log \Big(\frac K\rho \Big)\sum_{b \in B\colon
\mu(b)\alpha_t^d\le L}\mu(b)\Big\}\\
&\qquad\times \exp\Big\{\frac{t}{\alpha_t^{2}} \sum_{b \in B \colon \mu(b)\alpha_t^d
> L}\Big( M\frac K\rho\mu(b)-K \mu(b) \log\big( \mu(b)\alpha_t^d\big)\Big)\Big\}.
\end{aligned}
\end{equation}
Estimate the last term by
$$
\begin{aligned}
\frac{t}{\alpha_t^{d+2}} &\sum_{b \in B \colon \mu(b)\alpha_t^d > L}\biggl( M\frac
K\rho \mu(b)\alpha_t^d-K \mu(b)\alpha_t^d  \log\Big( \mu(b)\alpha_t^d\Big)\biggr)\\
&\leq\frac{t}{\alpha_t^{2}}\sum_{b \in B \colon \mu(b)\alpha_t^d >
L}\mu(b)\Big(M\frac K\rho-K\log L\Big),
\end{aligned}
$$
which is nonpositive for $L$ large enough (only depending on $M$, $K$ and $\rho$).
Now observe that $\sum_{b\in B}\mu(b)=1$ and that the cardinality of
$\mathcal{M}_1^{\ssup{D_t}}(B)$ is $\e^{o(t/\alpha_t^2)}$. This gives that
\begin{equation*}
\lim_{t\to\infty}\frac{\alpha_t^2}t\log\Big\langle F_{t,R}(\overline \xi_t\wedge M)^{K/\rho} \Big\rangle=\lim_{t\to\infty}\frac{\alpha_t^2}t\log\big(\mbox{l.h.s.~of }\eqref{lem-int1}\big)
\leq K\log \frac K\rho,
\end{equation*}
and the assertion (i) is proved.

Now we prove assertion (ii). We use the exponential Chebyshev inequality and (i) as
follows. We split the event $\{\overline\xi_t\wedge M\in D_{\beta,R}^{\rm
c}\}=\{|\Lcal_{3R}(\overline\xi_t\wedge M)-\rho|>\beta\}$ into the events
$\{\Lcal_{3R}(\overline\xi_t\wedge M)>\rho+\beta\}$. Let us consider only the first
of these events, the other is handled in the same way. On this event, we multiply
both sides of the inequality with $\e/\rho$, take logs, multiply with $\beta
t\alpha(t)^{-2}$ and take $\exp$. This gives, recalling the definition of
$F_{t,R}(\psi)$ in \eqref{Ftdef},
$$
\begin{aligned}
\1_{\{\Lcal_{3R}(\overline\xi_t\wedge M)>\rho+\beta\}}
&=\1_{\{F_{t,R}(\overline\xi_t\wedge M)^{\beta/\rho}>\exp\{t\alpha_t^{-2}\beta
\log(\e(1+\beta/\rho))\}\}}\\
&\leq F_{t,R}(\overline\xi_t\wedge M)^{\beta/\rho}\exp\Big\{-\frac
t{\alpha_t^{2}}\beta\log\big(\e\big(1+\sfrac\beta\rho\big)\big)\Big\}.
\end{aligned}
$$
Hence, we can estimate, with the help of assertion (i),
$$
\begin{aligned}
\big\langle F_{t,R}&(\overline\xi_t\wedge M) \1_{\{\Lcal_{3R}(\overline\xi_t\wedge
M)>\rho+\beta\}}\big\rangle\\
&\leq \big\langle  F_{t,R}(\overline\xi_t\wedge M)^{\beta/\rho+1}\big\rangle
\exp\Big\{-\frac
t{\alpha_t^{2}}\beta\log\big(\e\big(1+\sfrac\beta\rho\big)\big)\Big\}\\
&\leq \e^{o(t/\alpha_t^{2})} \exp\Big\{\frac
t{\alpha_t^{2}}\rho\Big[(\sfrac\beta\rho+1)\log(\sfrac\beta\rho+1)-\sfrac\beta\rho-\sfrac\beta\rho\log
\big(1+\sfrac\beta\rho\big)\Big]\Big\}\\
&=\e^{o(t/\alpha_t^{2})} \exp\Big\{\frac
t{\alpha_t^{2}}\rho\Big[\log(1+\sfrac\beta\rho) -\sfrac\beta\rho\Big]\Big\}.
\end{aligned}
$$
Since the term in square brackets is negative and  does not depend on $R$ nor on
$M$, the proof is complete.
\qed\end{Proof}

\section{The variational formulas}\label{sec-VarForm}

In this section we identify the constant $\chi$ appearing in Theorem~\ref{thm-HKM}
in terms of a \lq dual\rq\ variational problem which will be of importance.
Furthermore, we prove a minimisation property of that formula: every asymptotically
minimising sequence converges, along a suitable subsequence, after appropriate
spatial translation, towards the minimiser of the formula in $L^2(\R^d)$-sense. This
is one of the crucial ingredients of the subsequent proof of
Lemma~\ref{lem-apprvarform}, the last open step in the proof of
Proposition~\ref{ZweiteProp}.

Recall the parameter $\rho\in(0,\infty)$ from Assumption (HK).
Then \cite[Prop.~1.11]{HKM06} identifies $\chi$ as follows.

\begin{lemma}[Dual representation of $\chi$]\label{lemchirepr} For any $g\in H^1(\R^d)$, 
\begin{equation}\label{Hdef}
\Hcal(g^2)=\rho\int_{\R^d} g^2(x)\log g^2(x)\,\d x\;\in [-\infty,\infty)
\end{equation}
is well-defined. Furthermore, ${\mathcal L}$ and
$\mathcal{H}$ on $L^2(\mathbb{R}^d)$ are Legendre transform of each other, more
precisely,
\begin{equation}
\Lcal(\psi)=\sup_{g\in H^1(\R^d)}\big(\langle
g^2,\psi\rangle-\Hcal(g^2)\big)\qquad\mbox{and}\qquad
\Hcal(g^2)=\sup_{\psi\in\Ccal(\R^d)}\big(\langle g^2,\psi\rangle-\Lcal(\psi)\big).
\end{equation}
Furthermore,
\begin{equation}\label{Def_Chi_g}
\chi=\inf_{g \in H^1(\mathbb{R}^d)\colon ||g||_2=1} \Bigl\{\|\nabla
g\|_2^2-\mathcal{H}\bigl(g^2\bigl)\Bigr\}.
\end{equation}
Moreover, the minimum in \eqref{Def_Chi_g} is attained, uniquely up to translation,
at the Gaussian density
$$
\widehat g^2(x)=(\sfrac\rho\pi)^{\frac d2} \e^{-\rho |x|^2}=\sfrac 1\e
\e^{\frac 1\rho \widehat \psi(x)},\qquad x\in\R^d.
$$
The function $\widehat g$ is the unique $L^2$-normalized positive eigenfunction of the
operator $\Delta+\widehat\psi$ with eigenvalue $\lambda(\widehat\psi)=\rho-\rho
d+\rho\frac d2
\log\frac\rho\pi$ and satisfies $\Lcal(\widehat\psi)=\rho$.
\end{lemma}

The main point in the proof of Lemma~\ref{lemchirepr} is the well-known {\it
logarithmic Sobolev inequality},
\begin{equation}
 \label{logSobineq}
    \|\nabla g\|_2^2 \geq \Hcal(g^2)
    + \rho d\big(1-\sfrac 12 \log \sfrac \rho \pi\big),\qquad g\in L^2(\R^d),\|g\|_2=1,
\end{equation}
with equality if and only if $g$ is equal to $\widehat g$; see, e.g.,
\cite[Thm.~8.14]{LL01}. 

Now we consider the infimum in \eqref{Def_Chi_g} under the additional constraint
that any translation of $g^2$ is away from the minimizer $\widehat g^2$ introduced
in Lemma~\ref{lemchirepr} in $L^1(\R^d)$-sense, i.e.,
\begin{equation}\label{chiepsdef}
\chi(\eps)=\inf\Big\{\|\nabla g\|_2^2-\Hcal(g^2)\colon  g \in
H^1(\mathbb{R}^d),\|g\|_2=1, \forall x\in\R^d\colon \big\|{g}^2(x+\cdot)-\widehat
g^2(\cdot)\big\|_1\ge \eps \Big\}.
\end{equation}

The following lemma says that, given any $L^2$-normalised sequence $(g_n)_n$ of
approximative minimisers of $g\mapsto \|\nabla g\|_2^2-\Hcal(g^2)$, there is some
shift $x_n\in\R^d$ such that, along some subsequence, $g_n^2(x_n+\cdot)$ converges
in $L^1(\R^d)$ towards the Gaussian density $\widehat g^2$ introduced in
Lemma~\ref{lemchirepr}. Let us remark that a similar result is obtained in
\cite{C91} using a different approach. It is shown there that, for any
$L^2$-normalised $g\in H^1(\R^d)$,
\begin{equation}\label{Carlen}
\|\nabla g\|_2^2-2\pi \int g^2\log g^2\geq \chi +2\pi {\mathfrak
H}\big(|{\mathfrak F}g|^2\,\mathbf{\big|}\, 2^{d/2}\e^{-2\pi |x|^2}\big),
\end{equation}
where ${\mathfrak F}g(x)=\int_{\R^d} \e^{-2\pi \i\,x\cdot y}g(y)\,\d y$ is the
Fourier transform, and ${\mathfrak H}$ denotes the relative entropy between
probability measures with the respective densities. Note that the latter one is
equal to the Gaussian density $\widehat g^2$ with $\rho=2\pi$; by $\chi$ we mean 
our parameter with precisely that choice of $\rho$. Certainly,
\eqref{Carlen} can easily be generalised from $2\pi$ to any value of $\rho$.
However, the result in \eqref{Carlen} is not sufficient for our purposes since we found no way to
make the topology induced by this use of the entropy compatible with our
large-deviations arguments.

We also would like to mention that the discrete variant of the variational formula
in \eqref{Def_Chi_g} (i.e., where the Laplace operator in $\R^d$ is replaced by its
discrete version in $\Z^d$) has been analysed in detail in \cite[Theorem 2.II]{GH99}
and its dual variant in \cite[Proposition 1.1 and Lemma 3.2]{GKM05}. These are the
formulas that appear in the analysis of the parabolic Anderson model in the
universality class of the double-exponential distribution.

\begin{lemma}\label{lem-chiapprox} For any $\eps>0$, $\chi(\eps)>\chi$.
\end{lemma}
\begin{Proof}{Proof.}
It is sufficent to show that, for any sequence $(g_n)_n$ in $H^1(\R^d)$ such that
$\|g_n\|_2=1$ for all  $n$ and $\lim_{n\to\infty}(\|\nabla g_n\|_2^2-\rho \int
{g_n}^2 \log(g_n^2))= \chi$, there is a suitable shift $x_n\in\R^d$ such that, along
a suitable subsequence, $\lim_{n\to\infty}
\dist(g_n^2(x_n+\cdot),\widehat{g}^2(\cdot))=0$. Let $(g_n)_n$ be such a sequence.
Hence, for some $K>0$,
\begin{equation}
\|\nabla g_n\|_2^2-\rho \int {g_n}^2 \log({g_n}^2) \le K, \qquad  n \in \mathbb{N}.
\end{equation}
Now we show that $(\|\nabla g_n\|_2)_n$ is bounded: In the case {$d \ge 3$} we use
Jensen's inequality and the Sobolev inequality \cite[Theorem 8.3]{LL01} to estimate
\begin{equation}\label{energyesti}
\begin{aligned}
\|\nabla g_n\|^2_2 &\le K +\rho \int g_n^2 \log(g_n^2)
=  K +\rho \frac{d-2}{2} \int g_n^2 \log(g_n^{\frac{4}{d-2}})\\
&\le  K +\rho \frac{d-2}{2} \log\Big(\int  g_n^{\frac{2d}{d-2}}\Big)
\le  K +\rho \frac{d-2}{2} \log\Big(C \|\nabla  {g_n}\|^{\frac{2d}{d-2}}_2\Big),
\end{aligned}
\end{equation}
where $C>0$ is a Sobolev constant that satisfies $\int  f^{\frac{2d}{d-2}}\leq C \|\nabla  f\|^{\frac{2d}{d-2}}_2$ for any $f\in L^{\frac{2d}{d-2}}(\R^d)$. Hence $(\|\nabla g_n\|_2)_n$ is bounded and
therefore $(\int g_n^2 \log g_n^2)_n$ as well. In a similar way, we see the
boundedness of $(\|\nabla g_n\|_2)_n$ also in $d=2$, using the Sobolev inequality of
\cite[Theorem 8.5(ii)]{LL01}. In dimension $d=1$, we estimate, using the Sobolev
inequality of \cite[Theorem 8.5(i)]{LL01},
$$
\|\nabla g_n\|^2_2 \le K +\rho \int g_n^2 \log(g_n^2)\leq K+\rho\int
g_n^2\log\|g_n\|_\infty^2
\leq K+\rho\log\big(\sfrac 12\|\nabla g_n\|^2_2+\sfrac 12\big)
$$
and conclude as above.

Now we construct, for any $n\in\N$ and any small $\delta>0$ and any sufficiently
large $R=R_\delta>0$, some $x_n(\delta,R)\in\R^d$ such that 
\begin{equation}\label{gnyn}
\int_{x_n(\delta,R)+Q_R}g_n^2(y)\,\d y\geq 1-\delta.
\end{equation}
We pick a smooth auxiliary function $\Phi=\Phi_R\colon \R^d\to[0,1]$ satisfying
$\supp(\Phi)\subset Q_{R}$ and $\Phi\equiv 1$ on $Q_{R-1}$, and we put
$\Phi_x(y)=\Phi(x+y)$ for $x,y\in\R^d$. Consider $h_{n,x}=\Phi_x\cdot g_n$. Then we
have
$$
\int_{\R^d}\|h_{n,x}\|_2^2\,\d x=\int_{\R^d}\d y\,\int_{\R^d}\d
x\,\Phi^2(x+y)g_n^2(y)=\|\Phi\|_2^2.
$$
Similarly, we get
$$
\begin{aligned}
\int_{\R^d}\d x\,&\int_{\R^d}\d y\, h_{n,x}^2(y)\log (h_{n,x}^2(y))\\
&=\int_{\R^d}\d y\,g_n^2(y)\int_{\R^d}\d x\,\Phi^2(x+y)\log
(\Phi^2(x+y))+\|\Phi\|_2^2 \int g_{n}^2\log (g_{n}^2)\\
&=\int \Phi^2\log (\Phi^2)+\|\Phi\|_2^2 \int g_{n}^2\log (g_{n}^2).
\end{aligned}
$$
Using the product rule of differentiation, we get
$$
\begin{aligned}
\int_{\R^d}\|\nabla h_{n,x}\|_2^2\,\d x&=
\int_{\R^d}\d x\,\int_{\R^d}\d y\,\Big[g_n^2(y)|\nabla
\Phi(x+y)|^2+\Phi^2(x+y)|\nabla g_n(y)|^2
\\
&\qquad +2g_n(y) \Phi(x+y)\nabla \Phi(x+y)\cdot \nabla g_n(y)\Big]\\
&= \|\nabla \Phi\|_2^2 +\|\Phi\|_2^2 \|\nabla g_n\|_2^2+\int_{\R^d} g_n(y)\,
u_\Phi\cdot\nabla g_n(y)\,\d y,
\end{aligned}
$$
where $u_\Phi=\int_{\R^d}\d x\, \nabla(\Phi^2)(x)$. Using the
Cauchy-Schwarz inequality in the last term first for the Euclidean inner product and afterwards for the integral and recalling that $\|g_n\|_2=1$
and that $C=\sup_{n\in\N}\|\nabla g_n\|_2$ is finite, we see that
$$
\int_{\R^d}\|\nabla h_{n,x}\|_2^2\,\d x\leq \|\nabla \Phi\|_2^2 +\|\Phi\|_2^2
\|\nabla g_n\|_2^2+ C|u_\Phi|.
$$
Now fix $\delta>0$ and summarize the above estimates to obtain, for any $n\in\N$,
$$
\begin{aligned}
\int_{\R^d}\d x\,&\Big[ \|\nabla h_{n,x}\|_2^2-\rho \int_{\R^d}\d y\,
h_{n,x}^2(y)\log (h_{n,x}^2(y))
-(\chi+\delta)\|h_{n,x}\|_2^2\Big]\\
&\leq \|\nabla \Phi\|_2^2+\| \Phi\|_2^2\Big(\|\nabla g_n\|_2^2-\rho \int g_n^2\log
(g_n^2)-\chi-\delta\Big)+C|u_\Phi|-\rho \int\Phi^2\log(\Phi^2)\\
&=-\| \Phi\|_2^2(\delta-o(1))+\|\nabla \Phi\|_2^2+C|u_\Phi|-\rho
\int\Phi^2\log(\Phi^2),
\end{aligned}
$$
where $o(1)$ refers to $n\to\infty$ (recall that $(g_n)_n$ is asymptotically maximal
in the definition \eqref{Def_Chi_g} of $\chi$). It is possible to choose
$R=R_\delta$ so large that the right hand side is negative. Indeed, since $\Phi$
equals one in $Q_{R-1}$ and equals zero in $Q_R^{\rm c}$, all the terms $\|\nabla
\Phi\|_2^2$, $|u_\Phi|$ and $\int\Phi^2\log(\Phi^2)$ are of order $R^{d-1}$ as
$R\to\infty$, while $\| \Phi\|_2^2$ is of order $R^d$. Since the integral on the
left hand side is therefore also negative, there is some $x_n=x_n(\delta,R)\in\R^d$
such that the integrand is negative, that is,
$$
\Big\|\nabla\frac{h_{n,x_n}}{\|h_{n,x_n}\|_2}\Big\|_2^2-\rho\int
\Big(\frac{h_{n,x_n}}{\|h_{n,x_n}\|_2}\Big)^2\log\Big(\frac{h_{n,x_n}}{\|h_{n,x_n}\|_2}\Big)^2
\leq \chi+\delta+\rho\log\|h_{n,x_n}\|_2^2,
$$ 
using some elementary manipulations. By definition of $\chi$ in \eqref{Def_Chi_g},
the left hand side is no smaller than $\chi$, and it follows that
$\delta+\rho\log\|h_{n,x_n}\|_2^2$ is nonnegative. This in turn means that
$\|h_{n,x_n}\|_2^2\geq {\rm e}^{-\delta/\rho}\geq 1-\delta/\rho$. Replacing $
\delta/\rho$ by $\delta$, we have arrived at our first goal: the construction of
some $x_n(\delta,R)\in\R^d$ such that \eqref{gnyn} holds.

Now put $x_n=x_n(\frac 14,R_{1/4})$. We claim that the sequence
$(g_n^2(-x_n+\cdot))_{n\in\N}$ (conceived as probability measures on $\R^d$) is
tight. Indeed, for any $\delta\in(0,\frac 14)$ and any $n\in\N$, we pick $R_\delta$
and $x_n(\delta,R_\delta)$ as above. Since the masses of $g_n^2$ both in the box
$x_n(\frac 14,R_{1/4})+Q_{R_{1/4}}$ and in the box
$x_n(\delta,R_\delta)+Q_{R_\delta}$ exceed $\frac 34$, the two boxes must have an
non-empty intersection. Hence, the latter box is contained in the box $x_n(\frac
14,R_{1/4})+Q_{R_{1/4}+2R_\delta}$. Consequently, putting $\widetilde
R_\delta=R_{1/4}+2R_\delta$,
$$
\int_{Q_{\widetilde R_\delta}}g_n^2(-x_n+y)\,\d y
\geq \int_{x_n(\delta,R_\delta)+Q_{R_\delta}}g_n^2(y)\, \d y\geq 1-\delta.
$$
This shows the tightness of $(g_n^2(-x_n+\cdot))_{n\in\N}$. 

Now we use the Banach-Alaoglu theorem \cite[Theorem 2.18]{LL01} and \cite[Theorems
8.6, 8.7, 2.11]{LL01}. Since $(\|\nabla g_n(-x_n+\cdot)\|_2)_n$ is bounded, there is
a subsequence of $(g_n(-x_n+\cdot))_n$, still denoted   $(g_n(-x_n+\cdot))_n$, and a
${g} \in H^1(\R^d)$ satisfying $\|{g}\|_2\le 1$, such that $g_n(-x_n+\cdot)$
converges to ${g}$ weakly in $L^2(\R^d)$ and strongly in $L^p(Q_R)$ for any
$p<\frac{2d}{d-2}$ in $d\geq 3$ and for any $p<\infty$ in $d\in\{1,2\}$ and for any
$R\in(0,\infty)$ and almost everywhere, and $\nabla g_n(-x_n+\cdot)$ converges to
$\nabla {g}$ weakly in $L^2(\R^d)$. Furthermore, $\|\nabla g\|_2^2\le \liminf_{n \to
\infty}\|\nabla g_n\|_2^2$. Since $(g_n^2(-x_n+\cdot))_{n\in\N}$ is tight and by
local $L^2$-convergence, we also have that $g$ is $L^2$-normalized.

Now we argue that $\limsup_{n  \to \infty} \int g_n^2 \log(g_n^2)\leq \int {g}^2
\log{g}^2$. To derive this for $d\ge3$, we first estimate the integrals over
complements of large boxes. A similar estimate as the one in \eqref{energyesti}
shows, for any $R>0$ and $n\in\N$,
$$
\begin{aligned}
\int_{Q_R^{\rm c}} g_n^2(-x_n+y)\log (g_n^2(-x_n+y))\,\d y
&\leq -\int_{Q_R^{\rm c}} g_n^2(-x_n+y)\,\d y\log\Big(\int_{Q_R^{\rm c}}
g_n^2(-x_n+y)\,\d y\Big)\\
&\qquad +C\int_{Q_R^{\rm c}} g_n^2(-x_n+y)\,\d y,
\end{aligned}
$$
where $C>0$ is again a Sobolev constant. By tightness, the right hand side vanishes
as $R\to\infty$. A similar argument applies for $d\leq 2$.

Now we turn to the integral over the interior of a box. Observe that, for any $R>0$,
the sequence $g_n^2 \log g_n^2$ converges in probability to ${g}^2 \log{g}^2$ with
respect to the normalized Lebesgue measure on $Q_R$, that is, 
$$
\lim_{n \to \infty} \int_{Q_R} \d y
\,\1_{\{|g_n^2(-x_n+y)\log(g_n^2(-x_n+y))-{g}^2(y)\log({g}^2(y))|>\eta\}}=0,\qquad
\eta>0,
$$
as is easily deduced from the almost everywhere convergence of $g_n(-x_n+\cdot)$ to
${g}$, using Lebesgue's theorem.
Furthermore, $(g_n^2 \log(g_n^2))_n$ is uniformly integrable with respect to the
normalized Lebesgue measure on $Q_R$, which is seen, for $d \ge 3$, as follows. Note
that, for any $p\in(1 ,\frac{d}{d-2})$ and any $\beta\in(0,1)$, there is $c>0$ such
that 
\begin{equation}\label{xlogxesti}
|x \log(x)| \le c (|x|^p+ |x|^\beta), \qquad x >0, 
\end{equation}
and recall that $(\|g_n^2\|_{p'})_n$ and $(\|g_n\|_{2})_n$ are bounded for any $p'$
with  $p<p'<\frac{d}{d-2}$. From this it is easy to deduce the uniform integrability
on $Q_R$ for $d \ge 3$. A similar argument is used for $d=1,2$. Hence $\lim_{n  \to
\infty} \int_{Q_R} g_n^2 \log(g_n^2)= \int_{Q_R} {g}^2 \log({g}^2)$. Using the
above, we even see that $\limsup_{n  \to \infty} \int_{\R^d} g_n^2 \log(g_n^2)\leq 
\int_{\R^d} {g}^2 \log({g}^2)$.

Hence we see that $g$ is a minimizer in the definition \eqref{Def_Chi_g} of $\chi$.
Without loss of generality, we may therefore assume that $g$ is equal to $\widehat
g$ introduced in Lemma~\ref{lemchirepr}. Since $g_n(-x_n+\cdot)$ converges to
$\widehat g$ on every compact subset of $\R^d$ in $L^p$ for any $p\in(1,\frac
{2d}{d-2})$ in $d\geq 3$ and for any $p<\infty$ in $d\leq 2$, and by compactness of
$(g_n^2(-x_n+\cdot)$, we have also that $g_n^2(-x_n+\cdot)$ converges towards
$\widehat g$ in $L^2(\R^d)$-sense. This ends the proof.
\qed
\end{Proof}

Now we show that the variational formula $\chi(\eps)$ can be approximated by
finite-box versions. Introduce
\begin{equation}
\label{chiRepsdef}
\begin{aligned}
\chi_R(\eps)=\inf\Big\{\|\nabla g\|_2^2-\Hcal(g^2)\colon &g\in
H^1(\R^d),\|g\|_2=1,\supp(g)\subset Q_{3R},\\
&\forall x\in Q_{2R}\colon \int_{Q_R}|g^2(x+y)-\widehat g^2(y)|\,\d y\geq \eps\Big\}.
\end{aligned}
\end{equation}
Then $\chi_R=\chi_R(0)=\inf\{\|\nabla g\|_2^2-\Hcal(g^2)\colon g\in
H^1(\R^d),\|g\|_2=1,\supp(g)\subset Q_{3R}\}$.

\begin{lemma}[Finite-box approximation of $\chi$]\label{chiRapprox}
For any $\eps\geq 0$,
\begin{equation}
\liminf_{R\to\infty}\chi_R(\eps)\geq \chi(\eps).
\end{equation}
\end{lemma}

\begin{Proof}{Proof.} Let $(g_R)_{R\geq 1}$ be a family of $L^2$-normalised
functions $g_R\in H^1(\R^d)$ satisfying $\supp(g_R)\subset Q_{3R}$ and
$\int_{Q_R}|g_R^2(x+y)-\widehat g^2(y)|\,\d y\geq \eps$ for any $x\in Q_{2R}$ such
that $\|\nabla g_R\|_2^2-\Hcal(g_R^2)$ converges towards
$\liminf_{R\to\infty}\chi_R(\eps)$ as $R\to\infty$. Precisely as in the proof of
Lemma~\ref{lem-chiapprox}, we see that, for some sequence $R_n\to\infty$ as
$n\to\infty$, there are suitable shifts $x_n\in\R^d$ and some $L^2$-normalised $g\in
H^1(\R^d)$ such that $g_{R_n}(x_n+\cdot)$ converges towards $g$ in $L^2(\R^d)$ sense
and
\begin{equation}\label{chiRprop}
\liminf_{R\to\infty}\chi_R(\eps)=\lim_{n\to\infty}\big(\|\nabla
g_{R_n}\|_2^2-\Hcal(g_{R_n}^2)\big)
=\|\nabla g\|_2^2-\Hcal(g^2).
\end{equation}
By $L^2(\R^d)$-convergence of $g_{R_n}(x_n+\cdot)$ towards $g$, and since
$\int_{Q_R}|g_R^2(x+y)-\widehat g^2(y)|\,\d y\geq \eps$ for any $x\in Q_{2R}$, we
know that $g$ lies in the set of functions over which the infimum is taken in the
definition \eqref{chiepsdef} of $\chi(\eps)$. Hence, the right hand side of
\eqref{chiRprop} is not smaller than $\chi(\eps)$, which finishes the proof.
\end{Proof}
\qed

After these preparations, we finally can prove the last building block in the proof
of Proposition~\ref{ZweiteProp}.

\begin{Proof}{Proof of Lemma~\ref{lem-apprvarform}.} Recall the definition of
$\chi_R(\beta,\eps,t)$ from \eqref{chiapproxdef}; recall also
\eqref{rescaledeigenv}, \eqref{Dbetadef} and \eqref{GammaRepsdef}. We prove (i) and
(ii) jointly. Let $\psi_t\in \Fcal(Q_{3R})$, depending on $R>0$,
$\beta\in(0,\infty]$ and $\eps\geq 0$, be an approximatively maximizing function in
the definition \eqref{chiapproxdef} of $-\chi_R(\beta,\eps,t)$. More precisely, we
require that
$$
\lambda_{3R}^{\ssup t}(\psi_t)-\rho\log\Big( \frac\e\rho\Lcal_{3R}(\psi_t)\Big)\leq
-\chi_R(\beta,\eps,t)+\frac 1t.
$$ 
By the two extra conditions, we have 
$$
|\Lcal_{3R}(\psi_t)-\rho|\leq \beta\qquad\mbox{and}\qquad
\d_R(\psi_t(x+\cdot),\widehat\psi(\cdot))\geq \eps,\quad\mbox{for any }x\in Q_{2R}.
$$
By the first condition, we may pick some $c\in \R$ (to be precise,
$c=\rho\log(\rho/\Lcal_{3R}(\psi_t))$) such that 
\begin{equation}\label{psiass}
1=\frac 1\e\int_{Q_{3R}} \e^{\frac 1\rho[\psi_t(x)+c]}\,\d
x\qquad\mbox{and}\qquad|1-\e^{-c/\rho}|\leq \frac{\beta}\rho.
\end{equation}
Since $\lambda_{3R}^{\ssup t}(\psi_t+c)=\lambda_{3R}^{\ssup t}(\psi_t)+c$ and
$\Lcal_{3R}(\psi_t+c)=\e^{c/\rho}\Lcal_{3R}(\psi_t)$ and since 
$\Lcal_{3R}(\psi_t+c)=\rho$ by the choice of $c$, we have
$$
-\liminf_{t\to\infty}\chi_R(\beta,\eps,t)=\limsup_{t\to\infty}\lambda_{3R}^{\ssup
t}(\psi_t+c)-\rho.
$$
Recall the Rayleigh-Ritz formula $\lambdaD_R(V)=\max_{f\in\ell^2(B_R)\colon
\|f\|_2=1}(\langle \DeltaD f,f\rangle +\langle V,f^2\rangle)$ for potentials
$V\colon B_R\to\R$. Hence, there is an $\ell^2$-normalized function
$f_t\in\ell^2(\Z^d)$ in that is positive in $B_{3R\alpha_t}$ and zero outside and
satisfies
\begin{equation}
\lambda_{3R}^{\ssup t}(\psi_t+c)=\alpha_t^2
\lambdaD_{3R\alpha_t}\big(\sfrac{1}{\alpha_t^2}[\psi_t^{\rm d}+c]\big)
=\alpha_t^2 \langle\DeltaD {f}_t,{f}_t\rangle+\big\langle \psi_t^{\rm d}+c,
f_t^2\big\rangle. 
\end{equation}
For any $i\in\{1,\dots,d\}$, introduce $g_t^{\ssup i}\colon \R^d\to[0,\infty)$
defined by 
$$
g_t^{\ssup i}(x)=\alpha_t^{d/2}\Big[f_t(\lfloor x\alpha_t\rfloor)+\big(\alpha_t
x_i-\lfloor\alpha_t x_i\rfloor\big)\Big(f_t\big(\lfloor
x\alpha_t\rfloor+\e_i\big)-f_t\big(\lfloor x\alpha_t\rfloor\big)\Big)\Big],
$$ 
where $x=(x_i)_{i=1,\dots,d}$, and $\e_i\in\R^d$ is the $i$-th unit vector.
Abbreviate $\widetilde x_i=(x_j)_{j\not=i}\in\R^{d-1}$ and denote $g_{t,\widetilde
x_i}^{\ssup i}(x_i)=g_t^{\ssup i}(x)$. For almost every $\widetilde x_i\in\R^{d-1}$,
the map $ g_{t,\widetilde x_i}^{\ssup i}$ is continuous and piecewise affine, and
hence lies in $H^1(\R)$ with support in $[-3R,3R]$. Now let $(t_n)_{n\in\N}$ be a
sequence in $(0,\infty)$ with $\lim_{n\to\infty}t_n=\infty$ such that the limit
superior of $\lambda_{3R}^{\ssup t}(\psi_t+c)$ is realized along this sequence. Using Fubini's theorem and Fatou's lemma, one shows, in the same way as in the proof of \cite[Proposition ~5.1]{HKM06}, that
$$
\sum_{i=1}^d\int_{\R^{d-1}}d \widetilde x_i\,\liminf_{n\to\infty}\int_{\R} d x_i \,\big|(g^{\ssup{i}}_{n,\widetilde x_i})'(x_i)\big|^2<\infty.
$$
Furthermore, since $|x_i-\lfloor \alpha(t_n)x_i\rfloor/\alpha(t_n)|\leq \alpha(t_n)^{-1}$, one also derives that
$$
\lim_{n\to\infty}\|g_{t_n}^{\ssup{i}}-{\alpha(t_n)^{d/2} f_{t_n}(\lfloor\alpha(t_n)\,\cdot\,\rfloor)}\|_2=0.
$$
Hence, one sees that, along some subsequence, for almost every $\widetilde x_i\in\R^{d-1}$, $g_{t_n,\widetilde
x_i}^{\ssup i}$ converges towards some $g_{\widetilde x_i}^{\ssup i}\in H^1(\R^d)$.
According to \cite[Theorems~8.6 and 8.7]{LL01}, the convergence is strong in $L^q$ for any $q<\frac{2d}{d-2}$ for $d\geq 3$ and for
all $q<\infty$ for $d\in\{1,2\}$, pointwise almost everywhere and weak in $L^2$ for
the gradients. Furthermore, as also is shown in the proof of \cite[Proposition
~5.1]{HKM06}, there is some $L^2$-normalized $g\in H^1(\R^d)$ with support in
$Q_{3R}$ such that $g(x)=g_{\widetilde x_i}^{\ssup i}(x_i)$ for almost all
$x\in\R^d$, and we have
$$
\limsup_{t\to\infty}\alpha_t^2 \langle\DeltaD {f}_t,{f}_t\rangle\leq -\|\nabla g\|_2^2.
$$
Observe that, for any $i\in\{1,\dots,d\}$, 
$$
\begin{aligned}
\langle \psi_t^{\rm d}, f_t^2\rangle&=\int_{\R^d} (\psi_t(x)) f_t\big(\lfloor
\alpha_t x\rfloor\big)^2\alpha_t^d\,\d x\\
&=\int_{\R^d} (\psi_t(x))\Big(g_t^{\ssup i}(x)-\alpha_t^{d/2}\big(\alpha_t
x_i-\lfloor\alpha_t x_i\rfloor\big)\Big(f_t\big(\lfloor
x\alpha_t\rfloor+\e_i\big)-f_t\big(\lfloor x\alpha_t\rfloor\big)\Big)\Big)^2\,\d x.
\end{aligned}
$$
It is also clear from the proof of \cite[Proposition ~5.1]{HKM06} that the function
in the brackets on the right hand side has an $L^2$ distance to $g_t^{\ssup i}$ that vanishes
as $t\to\infty$ and that $g_t^{\ssup i}$ converges towards $g$ strongly in $L^2$. We
write now $g_t$ instead of $g_t^{\ssup 1}/\|g_t^{\ssup 1}\|_2$; recall that
$\lim_{t\to\infty}\|g_t^{\ssup 1}\|_2=1$. Hence, we have
$$
\begin{aligned}
-\liminf_{t\to\infty}\chi_R(\beta,\eps,t)&=\limsup_{t\to\infty}\lambda_{3R}^{\ssup
t}(\psi_t+c)-\rho\\
&\leq \limsup_{n\to\infty}\Big(\alpha_{t_n}^2 \langle\DeltaD
{f}_{t_n},{f}_{t_n}\rangle+\langle \psi_t^{\rm d}+c, f_{t_n }^2\rangle\Big)-\rho\\
&\leq -\|\nabla g\|_2^2+\limsup_{n\to\infty}\langle \psi_{t_n}+c-\rho,g_{t_n}^2\rangle.
\end{aligned}
$$
Now we employ the definition of $\Hcal$ in \eqref{Hdef} to rewrite
$$
\langle \psi_{t_n}+c-\rho,g_{t_n}^2\rangle= \Hcal(g_{t_n}^2)-\rho\Big\langle
g_{t_n}^2,\log\frac{g_{t_n}^2}{\e^{\frac 1\rho[\psi_{t_n}+c-\rho]}}\Big\rangle.
$$
Recall that $g_t$ is $L^2$-normalized and that $\e^{\frac1\rho[\psi_{t_n}+c-\rho]}$
is a probability density on $Q_{3R}$ by \eqref{psiass}. Hence, the last term is
equal to the entropy between the two probability measures with densities $\frac
1\e\e^{\frac1\rho[\psi_{t_n}+c]}$ resp.~$g_t^2$. According to
\cite[Ex.~6.2.17]{DZ98}, we can estimate this entropy against the variational
distance between these measures as follows.
$$
\Big\langle g_{t_n}^2,\log\frac{g_{t_n}^2}{\e^{\frac
1\rho[\psi_{t_n}+c-\rho]}}\Big\rangle
\geq \frac 12 \Big\|g_{t_n}^2-\frac1\e\e^{\frac1\rho[\psi_{t_n}+c]}\Big\|_{1,3R}^2,
$$
where $\|\cdot\|_{1,3R}$ denotes the $L^1$-norm on $L^1(Q_{3R})$. In the same way as
in the proof of Lemma~\ref{lem-chiapprox} (see around \eqref{xlogxesti}) one sees
that $\limsup_{n\to\infty}\Hcal(g_{t_n}^2)\leq \Hcal(g^2)$. Hence,
\begin{equation}\label{chiliminfbound}
-\liminf_{t\to\infty}\chi_R(\beta,\eps,t)\leq \Hcal(g^2)-\|\nabla
g\|_2^2-\frac\rho2\liminf_{n\to\infty}\Big\|g_{t_n}^2-\frac 1\e
\e^{\frac1\rho[\psi_{t_n}+c]}\Big\|_{1,3R}^2.
\end{equation}
Recall that $g\in H^1(\R^d)$ is $L^2$-normalized with support in $Q_{3R}$. Hence,
the right hand side can be estimated against $-\chi_R$, and this ends the proof of
Lemma~\ref{lem-apprvarform}(i), since we know from Lemma~\ref{chiRapprox} that
$\lim_{R\to\infty}\chi_R=\chi$. 

However, for proving (ii), we have to work harder in order to get an upper bound
that is strictly smaller. Recall the definition of $\chi_R(\eps)$ in
\eqref{chiRepsdef}. If $g$ is bounded away from $\widehat g$ in the sense that
$\int_{Q_R}|g^2(x+y)-\widehat g^2(y)|\,\d y\geq \frac \eps {4\e}$ for any $x\in
Q_{2R}$, then we can estimate the first two terms on the right hand side of
\eqref{chiliminfbound} from above against $-\chi_R(\frac\eps {4\e})$, which finishes
the proof of Lemma~\ref{lem-apprvarform}(ii), since
$\liminf_{R\to\infty}\chi_R(\frac \eps{4\e})\geq \chi(\frac \eps {4\e})>\chi$ by
Lemmas~\ref{chiRapprox} and \ref{lem-chiapprox}. Hence, it remains to consider the
case that $\int_{Q_R}|g^2(x+y)-\widehat g^2(y)|\,\d y< \frac \eps {4\e}$ for some
$x\in Q_{2R}$. Then we also have $\int_{Q_R}|g_{t_n}^2(x+y)-\widehat g^2(y)|\,\d y<
\frac \eps {2\e}$ for all sufficiently large $n$, since $g_{t_n}$ converges towards
$g$ in $L^2(Q_{3R})$. Now we estimate the last term on
  the right hand side of \eqref{chiliminfbound} as follows. Recall that $\widehat
g^2=\frac 1\e\e^{\widehat\psi/\rho}$ and use the reversed triangle inequality to
estimate 
\begin{equation}\label{reversedtriangle}
\begin{aligned} 
\Big\|g_{t_n}^2-\frac1\e\e^{\frac1\rho[\psi_{t_n}+c]}\Big\|_{1,3R}
&\ge \int_{Q_{R}}
\Big|\frac{1}{\e}\e^{\frac1\rho[(\psi_{t_n}(x+y))+c]}-g_{t_n}^2(x+y)\Big|\,\d y\\
&\ge \bigg|\e^{c/\rho}\int_{Q_{R}} \frac 1\e\Big|\e^{\frac1\rho \psi_{t_n}(x+y)}
-\e^{\frac1\rho\widehat\psi(y)}\Big|\,\d y\\
&\qquad-\big|\e^{c/\rho}-1\big|\int_{Q_R} \frac 1\e \e^{\frac
1\rho\widehat\psi(y)}\,\d y
-\int_{Q_{R}}\big|\widehat g^2(y)-g_{t_n}^2(x+y)\big|\,\d y\bigg|.
\end{aligned}
\end{equation}
Recall that $\psi_{t_n}\in \Gamma_{R,\eps}$, i.e., we have in particular that
$\int_{Q_R}\frac 1\e \big|\e^{\frac1\rho
\psi_{t_n}(x+y)}-\e^{\frac1\rho\widehat\psi(y)}\big|\,\d y\geq \frac\eps\e$, see \eqref{GammaRepsdef}.
Furthermore, we use the estimate for $c$ in \eqref{psiass} and the above mentioned
one for the distance between $g_{t_n}^2$  and $\widehat g^2$ to see that, for
$\beta>0$ small enough, the term between the outer absolute signs is positive and
may be estimated by 
\begin{equation}\label{lastesti}
\Big\|g_{t_n}^2-\frac1\e\e^{\frac1\rho[\psi_{t_n}+c]}\Big\|_{1,3R}
\geq \Big|\e^{c/\rho}\frac\eps\e -\e^{c/\rho}\frac\beta\rho-\frac\eps{2\e}\Big|.
\end{equation}
If one picks $\beta>0$ so small that $\e^{c/\rho}\geq 3/4$ and
$\e^{c/\rho}\frac\beta\rho\leq \frac\eps{8\e}$, then the right hand side of
\eqref{lastesti} is not bigger than $\frac\eps{8\e}$. This ends the proof.
\qed
\end{Proof}


\end{document}